\documentclass{article}
\usepackage{latexsym,amsmath,amssymb,amscd}
\usepackage[dvips]{color}

\begin{document}

\title{An Application of the Nash-Moser Theorem to 
the Vacuum Boundary Problem of Gaseous Stars }
\author{Tetu Makino \footnote{Professor Emeritus at Yamaguchi University, Japan / e-mail: makino@yamaguchi-u.ac.jp}}
\date{\today}
\maketitle

\newtheorem{Lemma}{Lemma}
\newtheorem{Proposition}{Proposition}
\newtheorem{Theorem}{Theorem}
\newtheorem{Definition}{Definition}

\begin{abstract}
We have been studying spherically symmetric motions of gaseous stars with physical
vacuum boundary governed either by the Euler-Poisson equations in the 
non-relativistic theory or by the Einstein-Euler equations in the relativistic theory. The problems are to construct solutions whose first approximations are small time-periodic solutions to the linearized problem at an equilibrium and to construct solutions to the Cauchy problem near an equilibrium. These problems can be solved when $1/(\gamma-1)$ is an integer, where $\gamma$ is the
adiabatic exponent of the gas near the vacuum, by the formulation by R. Hamilton of the Nash-Moser theorem. We discuss on an application of the formulation by J. T. Schwartz of the Nash-Moser theorem to the case in which 
$1/(\gamma-1)$ is not an integer but sufficiently large. \\


\end{abstract}

\section{Introduction}

In the previous works \cite{OJM}, \cite{KJM}, 
we investigated the time evolution of spherically symmetric gaseous stars,
either in the non-relativistic case governed by the Euler-Poisson
equations (\cite{OJM}), or in the relativistic
case governed by the Einstein-Euler equations (\cite{KJM}). Our studies suppose that the gas remains to be barotropic during the evolutions. That is,
the pressure $P$ is a given fixed function of the density $\rho$.
We were assuming \\

({\bf A}):\   {\it  $P$ is a given smooth function of
$\rho>0$ such that $P>0, dP/d\rho>0$ for $\rho>0$ and there are positive constants $A, \gamma$ such that $1<\gamma<2$ and an analytic function
$\Omega$ on a neighborhood of $0$ such that $\Omega(0)=1$ and
$$P=A\rho^{\gamma}\Omega(\rho^{\gamma-1})$$
for $0<\rho\ll 1$. }\\

Assuming that there is an equilibirium
with a finite radius at which the gas touches the vacuum, we investigated time-dependent solutions 
near this equilibirium. The existence of solutions whose first approximations are
small time-periodic solutions to the linearized problem at the equilibrium, and the
existence of solutions to the Cauchy problem near this equilibrium
were established by applying the Nash-Moser theorem formulated by
R. Hamilton, \cite{Hamilton}, which we shall call the {\bf `Nash-Moser(-Hamilton) theorem'}. 

But in order to apply this Nash-Moser(-Hamilton) theorem, we had to put the assumption\\

({\bf B}):\  {\it $N$ is an even integer.}\\

\noindent Here the parameter $N$ is determined from the approximate
adiabatic exponent $\gamma$ near the vacuum by
$$\frac{N}{2}=1+\frac{1}{\gamma-1},\qquad\mbox{or}
\qquad \gamma=1+\frac{2}{N-2}. $$
Under the assumption ({\bf B}), the function $(1-x)^{N/2}$ is
analytic at $x=1$, and the smoothness
of this function plays an essential r\^{o}le
for the application of the Nash-Moser(-Hamilton) theorem.

However in many physically important cases, in which, e.g., 
$\gamma=5/3, 7/5$ and so on, $N/2$ is not an integer. Therefore
the open problem to apply the Nash-Moser theorem for the case
in which $N$ is not an even integer is very important. 

The present study is a partial answer to this open problem. In fact,
when $N$ is very large, that is, $\gamma-1$ is very small,
the Nash-Moser theorem formulated by J. T. Schwartz, \cite{Schwartz}, 
which we shall call the {\bf `Nash-Moser(-Schwartz) theorem'}, can be applied. To show this is the aim of this article.\\

 The reason why we apply the Nash-Moser theorem is that we have to treat the so called `physical vacuum boundary', that is, a boundary 
$\partial\Omega$ of the domain $\Omega$ on which $\rho>0$ and 
outside which $\rho=0$ such that 
$$ 0<-\frac{\partial}{\partial\mathbf{N}} 
\Big(\frac{dP}{d\rho}\Big)<+\infty,$$
where $\mathbf{N}$ is the outer normal vector of $\partial\Omega$. 
(Cf. \cite{CoutandS2012}, \cite{JangM2009} ). Because of this singularity at the physical vacuum boundary, the nonlinear hyperbolic
evolution equation to be considered involves a loss of the derivative
regularities by which the usual iteration does not work. (Cf.
\cite[p. 49]{FE}.) This difficulty has already been attacked by several scholars: D. Coutand and S. Shkoller, \cite{CoutandS2011}, \cite{CoutandS2012}, J. Jang and N. Masmoudi, \cite{JangM2009},
\cite{JangM2015}, and T. Luo, Z.-P. Xin and H.-H. Zen \cite{LuoXZ2014}.
Although they do not consider the approximation by time-periodic
solutions of the linearized problem, they have already developed 
a powerful theoretical method based on sophisticated energy estimates
using the Hardy inequality. Even for the time-periodic approximation, 
a work along this line has been done by J. Jang, \cite{Jang2016}. 
Maybe their method, after a suitable generalization, 
can be applied to the Einstein-Euler system, too. The advantage of the 
method developed by them is the absence of additional restrictions on
the number theoretic properties or the nearness to $1$ of the adiabatic exponent $\gamma$,
which seems inevitable if we insist on the application of the 
Nash-Moser technique.

Thus, in this situation, the present study may give, as an another approach, a touch stone for their method.

\section{Problem setting}

We are going to consider the following problem.\\

The equations to be considered are 
\begin{subequations}
\begin{align}
&\frac{\partial y}{\partial t}-J\Big(x, y, x\frac{\partial y}{\partial x}\Big)v=0, \label{G01a}\\
&\frac{\partial v}{\partial t}+H_1\Big(x, y, x\frac{\partial y}{\partial x},
 v\Big)\mathcal{L}\Big(\frac{\partial}{\partial x}\Big)y+
H_2\Big(x,y, x\frac{\partial y}{\partial x}, v, x\frac{\partial v}{\partial x}
\Big)=0 \label{G01b}
\end{align}
\end{subequations}
with
\begin{equation}
\mathcal{L}\Big(\frac{d}{dx}\Big)=-x(1-x)\frac{d^2}{dx^2}-
\Big(\frac{5}{2}(1-x)-\frac{N}{2}x\Big)\frac{d}{dx}+
L_1(x)\frac{d}{dx}+L_0(x). \label{G02}
\end{equation}

The independent variables $t, x$ run over $0\leq t\leq T, 0<x<1$, where
$T$ is a fixed number. The unknown functions are $y(t,x), v(t,x)$.\\

 In fact, in the case of the relativistic Einstein-Euler system
discussed in \cite{KJM}, the equations (\ref{G01a})(\ref{G01b})(\ref{G02})
are derived as follows.
The metric to be fonund is
$$ds^2
=e^{2F}c^2dt^2
-e^{2H}dr^2
-R^2(d\theta^2+\sin^2\theta d\phi^2), $$
and the unknown density distribution is $\rho=\rho(t,r)$. We fix an equilibrium governed by the Tolman-Oppenheimer-Volkoff equation
$\rho=\bar{\rho}(r)$ with
$$\bar{m}(r):=4\pi \int_0^r\bar{\rho}(r')r'^2dr',$$
for which the radius $r_+$ and the total mass $m_+$ are finite,
and we introduce the perturbation variables $y, v$ by
$$R=r(\bar{m})(1+y),\qquad
V=r(\bar{m})v,$$
by taking $\bar{m}(r)$ as the independent variable. Here
$$V=e^{-F}\Big(\frac{\partial R}{\partial t}\Big)_{r=\mbox{Const.}}.$$
By a suitable change of the variable from $\bar{m}$ (running on the interval $[0,m_+]$) to $x$ (running on $[0,1]$), the equations to be considered are reduced to (\ref{G01a})(\ref{G01b})(\ref{G02}). 
Actually we have
\begin{align*}
\rho&=\bar{\rho} (1+y)^{-2}\Big(1+y+
\bar{r}\frac{\partial y}{\partial\bar{r}}\Big)^{-1}, \\
e^F&=\sqrt{\kappa}\exp\Big[-\frac{u}{c^2}\Big], \\
J&=e^F(1+P/c^2\rho),
\end{align*}
and so on. Here $\kappa:=1-2Gm_+/c^2r_+$ and
$$u:=\int_0^{\rho}\frac{dP}{\rho+P/c^2}.$$
For the details, see \cite{KJM}

 In the case of the non-relativistic Euler-Poisson system, we can
 take $J=1$, naturally since $c=+\infty$, and the system of the
 first order equations (\ref{G01a})(\ref{G01b}) is reduced to
 a single second order equation of the form
 $$\frac{\partial^2y}{\partial t^2} +H_1\mathcal{L}y+H_2=0.$$
 For the details, see \cite{OJM}.\\
 
 Anyway let us list up the assumptions on the
 system (\ref{G01a})(\ref{G01b})(\ref{G02}),
 which have been ascertained for the Euler-Poisson system
 in \cite{OJM} and for the Einstein-Euler system in \cite{KJM}.\\

First we assume \\

 {\bf (B0):} $N\geq 5$.\\

Let us denote by $\mathfrak{A}_{(N)}$ the set of
all smooth functions $f(x)$ of $x \in [0,1[$ such that
$$
f(x)=
\begin{cases}
[x]_0 \quad\mbox{as}\quad x\rightarrow 0 \\
[1-x,(1-x)^{N/2}]_0 \quad\mbox{as}\quad x\rightarrow 1.
\end{cases}
$$

Here and hereafter $[X]_Q$ stands for a convergent power series of the
form $\displaystyle \sum_{j\geq Q}a_jX^j$, and
$[Z_1,Z_2]_Q$ stands for a convergent double power series of the form
$\displaystyle \sum_{j_1+j_2\geq Q}a_{j_1j_2}Z_1^{j_1}Z_2^{j_2}$.

By $\mathfrak{A}_{(N)}^Q(U^p)$,
$U$ being a neighborhood of $0$, we denote the set of
all smooth functions 
$f(x,y_1,\cdots, y_p)$ of
$x \in [0,1[, y_1,\cdots. y_p \in U$ such that 
there are convergent power series
$$\Phi_0(X,Y_1,\cdots, Y_p)=\sum_{k_1+\cdots+ k_p\geq Q}
a^{[0]}_{jk_1\cdots k_p}X^jY_1^{k_1}\cdots Y_p^{k_p}$$
and $$\Phi_1(Z_1, Z_2, Y_1,\cdots, Y_p)=
\sum_{k_1+\cdots+k_p\geq Q}
a^{[1]}_{j_1j_2k_1\cdots k_p}
Z_1^{j_1}Z_2^{j_2}Y_1^{k_1}\cdots Y_p^{k_p}$$
such that
$$f(x,y_1,\cdots, y_p)=\Phi_0(x,y_1,\cdots, y_p) \quad\mbox{for}\quad 0<x\ll 1$$
and 
$$
f(x, y_1,\cdots, y_p)=\Phi_1(1-x, (1-x)^{N/2},y_1,\cdots, y_p)
\quad\mbox{for}\quad 0<1-x\ll 1.$$
Using these notations we assume \\

{\bf (B1; $N$):} {\it $L_0, L_1\in \mathfrak{A}_{(N)} $ and
$$
L_1(x)=
\begin{cases}
[x]_1\quad\mbox{as}\quad x\rightarrow 0 \\
[1-x,(1-x)^{N/2}]_1\quad\mbox{as}\quad x\rightarrow 1
\end{cases}
$$
and
there is a neighborhood $U$ of $0$ such that
$J\in \mathfrak{A}_{(N)}^0(U^2)$,
$H_1\in\mathfrak{A}_{(N)}^0(U^3)$,
$H_2\in\mathfrak{A}_{(N)}^2(U^4)$.}\\

\textbullet \hspace{5mm} Note that, if $f(x)$ is a function of $\mathfrak{A}_{(N)}$, there is a smooth function
$\Phi \in C^{\infty}([0,1]\times[0,1])$ such that
$$f(x)=\Phi(1-x, (1-x)^{N/2})\qquad \mbox{for}\quad 0\leq\forall x <1, $$
althogh such an analytic two variable function $\Phi$
may not exist.\\

We suppose the following assumptions {\bf (B2)}, {\bf (B3; $N$)}:\\

{\bf (B2):}\  We have $J(x, 0, 0)H_1(x, 0, 0, 0)=1$ and there is a constant $C$ 
such that 
$$\frac{1}{C}<J(x,0,0)<C.$$\\

{\bf (B3; $N$):}\  We have
$$\partial_zJ\equiv_{(N)}0,\quad
(\partial_zH_1)\mathcal{L}y+\partial_zH_2\equiv_{(N)}0,
\quad \partial_wH_2\equiv_{(N)}0 $$
as $x\rightarrow 1$. \\

Here $z, w$ stand for $\displaystyle x\frac{\partial y}{\partial x},
x\frac{\partial v}{\partial x} $ respectively, and
 ``$f\equiv_{(N)} 0$ as $x\rightarrow 1$" is defined as that there is 
a convergent power series $\Phi(Z_1,Z_2, Y_1,\cdots, Y_p)$ such that
$$f(x,y_1,\cdots, y_p)=(1-x)
\Phi(1-x, (1-x)^{N/2}, y_1, \cdots, y_p)
\quad\mbox{for}\quad 0<1-x\ll 1.
$$\\

We assume \\

$\neg${\bf (B)} : $N$ is not an even integer. \\

So we consider the equations (\ref{G01a})(\ref{G01b}) under the assumptions
{\bf (B0), (B1; $N$), (B2), (B3; $N$)} and {\bf $\neg$(B)}.\\

Note that the result of the spectral analysis of the linear operator $\mathcal{L}$ is the same as the case in which $N/2$ is
an integer. See Appendix.\\

Hence the problem may be settled as follows:\\

Fix $T>0$, and functions $y^*, v^* \in \mathfrak{B}_{(N)}([0,T]; U)$.
Here $\mathfrak{B}_{(N)}([0,T]; U)$ stands for the set of all smooth functions $u(t,x)$ of $0\leq t\leq T, 0\leq x<1$, valued in $U$, such that
there are analytic functions $\Phi_0$
on $[0,T]\times]-\delta,\delta[$ and $\Phi_1$ 
on $[0,T]\times ]-\delta,\delta[^2$ with $0<\delta\ll 1$ such that
$$ u(t,x) =\Phi_0(t, x) \qquad\mbox{for}\quad 0\leq t\leq T, 0<x\ll 1$$ and
$$ u(t,x)=\Phi_1(t, 1-x, (1-x)^{N/2}) \quad\mbox{for}\quad 0\leq t\leq T, 0<1-x\ll 1$$. \\ 

Then we seek a solution $(y,v)$ of (\ref{G01a})(\ref{G01b}) of the form
\begin{equation}
y=y^*+\tilde{y}, \qquad v=v^*+\tilde{v} \label{G03}
\end{equation}
such that
\begin{equation}
\tilde{y}|_{t=0}=0, \qquad \tilde{v}|_{t=0}=0. \label{G04}
\end{equation}

Actually in the application to the problem of spherically symmetric gaseous stars, we consider the following $(y^*, v^*)$:\\

{\bf [1): To construct solutions whose first approximations are small time-periodic
solutions of the linearized problem]} Let $\lambda$ be a positive eigenvalue of $\mathcal{L}$ and $\Phi(x)$ be
the associated eigenfunction (see Appendix); put
\begin{align*}
Y_1(t,x)&=\sin(\sqrt{\lambda}t+\Theta_0)\Phi(x), \\
V_1(t,x)&=
\frac{\sqrt{\lambda}}{J(x,0,0)}\cos(\sqrt{\lambda}t+\Theta_0)\Phi(x),
\end{align*}
$\Theta_0$ being a constant; 
let $\varepsilon$ be small parameter; put $y^*=\varepsilon Y_1,
v^*=\varepsilon V_1$.\\

{\bf [2): To construct solutions of the Cauchy problem]} 
Let $\psi_0(x)=y|_{t=0}, \psi_1(x)=v|_{t=0}$ be the smooth initial data;
put
$$y^*=\psi_0(x)+tJ(x,0,0)\psi_1(x),\quad v^*=\psi_1(x).$$\\

{\bf Remark.} \hspace{7mm} In the study of the Einstein-Euler equations, \cite{KJM},
we fixed a short equilibrium $\rho(r), 0<r<r_+,$ for which we showed that 
$$u(r):=\int_0^{\rho(r)}\frac{dP}{\rho+P/c^2} $$
is analytic in $r$ at $r=r_+$ under the assumption ({\bf B}). (See \cite[Proposition 6]{KJM}.) However
we can prove that, even if $\neg$({\bf B}) is supposed,
we have
$$u(r)=B(r_+-r)(1+[r_+-r, (r_+-r)^{N/2}]_1) $$
as $r\rightarrow r_+$ with a positive constant $B$
for any short equilibrium. Actually a proof can be 
found in \cite{TOVdS}. (See \cite[Theorem 4]{TOVdS}.) But at the
moment, we may assume this property for the fixed equilibrium
under consideration. This analytic property of the equilibrium at the
vacuum boundary guarantees the validity 
of ({\bf B1}; $N$) and ({\bf B3}; $N$) for our application to the
study of spherically symmetric gaseous stars.

\section{Nash-Moser(-Schwartz) theorem}

We are going to apply the following Nash-Moser(-Schwartz) theorem 
(See \cite[Chapter II]{Schwartz}
and also see \cite{Sergeraert}):\\

{\bf Let $J$ be a positive integer and
$E_j, j=0,1,\cdots, J, F_j, j=1, \cdots, J$ be Banach spaces
such that $E_{j+1}\subset E_j, \|\mathfrak{u}\|_{E_j}\leq
\|\mathfrak{u}\|_{E_{j+1}},$ $F_{j+1}\subset F_{j},
\|\mathfrak{u}\|_{F_j}\leq \|\mathfrak{u}\|_{F_{j+1}}$.

For $\theta \geq 0$ there is a linear operator $S(\theta): E_0\rightarrow E_J$,
so called a smoothing operator, such that, if $0\leq j\leq \bar{j}\leq J$, then
\begin{align*}
&\|S(\theta)\mathfrak{u}\|_{E_{\bar{j}}}\leq C\theta^{\bar{j}-j}\|\mathfrak{u}\|_{E_j}, \\
&\|(I-S(\theta))\mathfrak{u}\|_{E_j}\leq C
\theta^{j-\bar{j}}\|\mathfrak{u}\|_{E_{\bar{j}}}.
\end{align*}

$\mathcal{F}$ is a mapping from $V:=\{\mathfrak{u}\in E_1\  |\  \|\mathfrak{u}\|_{E_1}<1\}$ into
$F_1$ such that, for $1\leq j \leq J$, $\mathcal{F}(V\cap E_j)\subset F_j$.
For any $\mathfrak{u}\in V$ there is a linear operator $D\mathcal{F}(\mathfrak{u}): E_1\rightarrow F_1$ such that
\begin{align*}
&\|D\mathcal{F}(\mathfrak{u})\mathfrak{h}\|_{F_1}\leq C
\|\mathfrak{h}\|_{E_1}\quad\mbox{for}\quad \mathfrak{u}\in V, \mathfrak{h}\in E_1, \\
&\|\mathcal{F}(\mathfrak{u}+\mathfrak{h})-\mathcal{F}(\mathfrak{u})-D\mathcal{F}(\mathfrak{u})\mathfrak{h}\|_{F_1}\leq
C\|\mathfrak{h}\|_{E_1}^2
\quad\mbox{for}\quad \mathfrak{u}, \mathfrak{u}+\mathfrak{h}\in V.
\end{align*}

For $\mathfrak{u}\in V$ there is a linear operator $\mathfrak{I}(\mathfrak{u}):F_1\rightarrow E_0$ such that, for $\mathfrak{u}\in V\cap E_j, 1\leq j\leq J$, $\mathfrak{I}(\mathfrak{u})F_j
\subset E_{j-1}$ and
\begin{align*}
&\|\mathfrak{I}(\mathfrak{u})\mathfrak{g}\|_{E_0}\leq C\|\mathfrak{g}\|_{F_1}\quad\mbox{for}\quad \mathfrak{u}\in V, \mathfrak{g}\in F_1, \\
&D\mathcal{F}(\mathfrak{u})\mathfrak{I}(\mathfrak{u})\mathfrak{g}=\mathfrak{g}\quad\mbox{for}\quad
\mathfrak{u}\in V, \mathfrak{g}\in F_J, \\
&\|\mathfrak{I}(\mathfrak{u})\mathcal{F}(\mathfrak{u})\|_{E_{J-1}}\leq C(1+\|\mathfrak{u}\|_J)\quad\mbox{for}\quad \mathfrak{u} \in V\cap E_J.
\end{align*}

Then there is a small $\delta$ such that if
$\|\mathcal{F}(0)\|_{E_1}\leq\delta$ then there is a solution $\mathfrak{u}$
of the equation $\mathcal{F}(\mathfrak{u})=0$ in $V$, provided that $J\geq 10$. } \\

We are going to apply this theorem to the equation
$\mathfrak{F}(\vec{w})=0$, where
$\mathfrak{F}(\vec{w})=(F_1, F_2)^T$ is the left-hand side of
(\ref{G01a})(\ref{G01b}) with $y=y^*+\tilde{y}, v=v^*+\tilde{v}$,
$\vec{w}=(\tilde{y},\tilde{v})^T$.\\

We will take the graded Banach spaces of functions
$\vec{u}(t,x)$ of $0\leq t\leq T, 0\leq x\leq 1$ defined by the norms
$$\|\vec{u}\|_{E_j}=\|\vec{u}\|_{\mathfrak{b}_E+\mathfrak{r}j}^{(2)},
\qquad
\|\vec{u}\|_{F_j}=\|\vec{u}\|_{\mathfrak{b}_F+\mathfrak{r}j}^{(2)}.$$
Later the positive parameters $\mathfrak{b}_E, \mathfrak{b}_F, \mathfrak{r}$ will be chosen suitably.
Here the norms $(\|\cdot\|_{\nu}^{(2)})_{\nu}$ are defined as follows:\\

 For a function $u(x)$ of $0\leq x\leq 1$, we put
$$u^{[0]}(x)=\omega(x)u(x),\qquad
u^{[1]}(x)=(1-\omega(x))u(x), $$
where
$\omega \in C^{\infty}(\mathbb{R})$ such that
$\omega(x)=1$ for $x\leq 1/3$, $0<\omega(x)<1$ for
$1/3<x<2/3$, and
$\omega(x)=0$ for $2/3\leq x$. 

Let us denote
\begin{align*}
&\triangle_{[0]}=x\frac{d^2}{dx^2}+\frac{5}{2}\frac{d}{dx},
\quad \triangle_{[1]}=X\frac{d^2}{dX^2}+\frac{N}{2}\frac{d}{dX}\quad\mbox{with}\quad X=1-x, \\
& \dot{D}_{[0]}=\sqrt{x}\frac{d}{dx},\quad \dot{D}_{[1]}=\sqrt{X}\frac{d}{dX}\quad
\mbox{with}\quad X=1-x, \\
&\|u\|_{[0]}=\Big(\int_0^1|u(x)|^2x^{3/2}dx\Big)^{1/2}, \\
&\|u\|_{[1]}=\Big(\int_0^1|u(x)|^2X^{N/2-1}dX\Big)^{1/2}.
\end{align*}

We put, for $\mu=0,1$, 
\begin{align*}
&\langle u\rangle _{[\mu]\ell}=
\begin{cases}
\|\triangle_{[\mu]}^mu\|_{[\mu]}\quad\mbox{for}\quad \ell=2m, \\
\|\dot{D}_{[\mu]}\triangle_{[\mu]}^mu\|_{[\mu]}\quad\mbox{for}\quad
\ell=2m+1,
\end{cases}\\
&\|u\|_{[\mu]k}=\Big(
\sum_{0\leq\ell\leq k}(\langle u\rangle_{[\mu]\ell})^2\Big)^{1/2}, \\
&\|\vec{u}\|_{[\mu]k}=\Big((\|y\|_{[\mu]k+1})^2+
(\|v\|_{[\mu]k})^2\Big)^{1/2} 
\end{align*}
for $\vec{u}=(y,v)^T=(y(x), v(x))^T$,
\begin{align*}
&\|u\|_k=\Big( (\|u^{[0]}\|_{[0]k})^2+(\|u^{[1]}\|_{[1]k})^2
\Big)^{1/2}, \\
&\|\vec{u}\|_k=\Big( (\|y\|_{k+1})^2+(\|v\|_k)^2
\Big)^{1/2}
\end{align*}
for $\vec{u}=(y,v)^T$.

For a function $\vec{u}=\vec{u}(t,x)$ of
$0\leq t\leq T, 0\leq x\leq 1$, we put
\begin{align*}
&|\vec{u};\tau, n\|_{[\mu]}      =\sup_{0\leq t\leq \tau}
\sum_{j+k\leq n}
\|\partial_t^j\vec{u}(t,\cdot)\|_{[\mu]k}, \\
&|\vec{u};\tau,n\| =\Big( (|\vec{u}^{[0]}; \tau,n\|_{[0]})^2+
(|\vec{u}^{[1]}; \tau,n\|_{[1]})^2
\Big)^{1/2},\\
&|\|\vec{u}\||_{[\mu]n}=\Big(
\sum_{j+k\leq n}\int_0^T(\|\partial_t^j\vec{u}(t,\cdot)\|_{[\mu]k})^2dt
\Big)^{1/2}, \\
&|\|\vec{u}\||_n=\Big(
(|\|\vec{u}^{[0]}\||_{[0]n})^2+(|\|\vec{u}^{[1]}\||_{[1]n})^2
\Big)^{1/2}.
\end{align*}

On the other hand, we put
\begin{align*}
&\|u\|_{[\mu]\kappa}^*=\Big(
\sum_{0\leq m\leq\kappa}(\|\triangle_{[\mu]}^mu\|_{[\mu]})^2
\Big)^{1/2}, \\
&\|u\|_{\kappa}^*=\Big(
(\|u^{[0]}\|_{[0]\kappa}^*)^2+(\|u^{[1]}\|_{[1]\kappa}^*)^2
\Big)^{1/2}, \\
&\|u\|_{\nu}^{(2)}=\Big(
\sum_{\iota+\kappa\leq\nu}\int_0^T
(\|(-\partial_t^2)^{\iota}u(t,\cdot)\|_{\kappa}^*)^2dt
\Big)^{1/2}, \\
&\|\vec{u}\|_{\nu}^{(2)}=\Big(
(\|y\|_{\nu}^{(2)})^2+(\|v\|_{\nu}^{(2)})^2
\Big)^{1/2}
\end{align*}
for $\vec{u}=(y,v)^T$.\\

{\bf Notation:}\    Hereafter we shall denote
$X\lesssim Y$ for two quantities $X, Y$ if there is a constant
$C$ such that $X \leq CY$. We shall denote $X\simeq Y$ if both
$X \lesssim Y$ and $ X\gtrsim Y$ hold, that is,
there is a constant $C$ such that
$(1/C)Y\leq X\leq CY$.\\

We can claim
\begin{equation} |\|\vec{u}\||_n\lesssim |\vec{u}; T,n\|
\lesssim |\|\vec{u}\||_{n+1}, \label{X01}
\end{equation}
and
\begin{equation}
|\|\vec{u}\||_{2\nu-1}\lesssim\|\vec{u}\|_{\nu}^{(2)}\lesssim
 |\|\vec{u}\||_{2\nu}. \label{X02}
\end{equation}

In fact we have
\begin{align*}
&\|\dot{D}_{[\mu]}u\|_{[\mu]}\lesssim
\|u\|_{[\mu]}+\|\triangle_{[\mu]}u\|_{[\mu]}, \\
&\int_0^T(\partial_t u)^2dt\lesssim
\int_0^Tu^2dt+
\int_0^T(-\partial_t^2u)^2dt
\end{align*}
and so on.

\section{Smoothing operators}

We construct the smoothing operator
\begin{align*}
&\vec{S}(\theta)\vec{u}=(S(\theta)y, S(\theta)v)^T
\quad\mbox{for}\quad \vec{u}=(y,v), \\
&S(\theta)=S_{[0]}(\theta)u^{[0]}+
S_{[1]}(\theta)u^{[1]}
\end{align*}
as follows.

Let us define $S_{[1]}(\theta)$.
For the simplicity, we write $S=S_{[1]}, \triangle=\triangle_{[1]},
$ $\|\cdot\|=\|\cdot\|_{[1]}, x$ instead of $1-x$, so that,
$\triangle =xD^2+(N/2)D, D=d/dx$, and
$u(x)$ be a function which vanishes for $2/3\leq x$.

We fix an extension of functions $u\in C^{\infty}([0,T]\times[0, 1[)$ such
that $u(t,x)=0$ for $2/3\leq x$ to a function
$\tilde{u}\in C_0^{\infty}(]-2T, 2T[\times[0, 1[)$ such that
$$\|\tilde{u}\|_{\nu}^{\sharp}\simeq \|u\|_{\nu}^{(2)},$$
where
$$\|\tilde{u}\|_{\nu}^{\sharp}:=
\Big(\sum_{\iota+\kappa\leq\nu}
\int_{-2T}^{2T}(\|(-\partial_t^2)^{\iota}\triangle^{\kappa}\tilde{u}(t,\cdot)
\|)^2dt
\Big)^{1/2}.
$$
See \cite[\S 5.17-30]{AdamsF}.

Let $(\phi_a(t))_{a=1,2,\cdots}$ be the eigenfunctions of the operator
$-\partial_t^2$ in $L^2(]-2T, 2T[)$ with the Dirichlet boundary conditions at $t=\pm 2T$ and $(\lambda_a)_a$ be the associated eigenvalues. That is,
$$\phi_a(t)=\frac{1}{\sqrt{2T}}\sin\frac{a\pi t}{2T},
\qquad
\lambda_a=\Big(\frac{a\pi}{2T}\Big)^2. $$
Hence $\lambda_a \simeq a^2$.

Let $(\psi_b(x))_{b=1,2,\cdots}$ be the eigenfunctions of $-\triangle$
in $L^2(x^{N/2-1}dx)$ with the Dirichlet boundary condition
at $x=1$ and $(\mu_b)_b$ be the associated eigenvalues.
That is,
$\psi_b(x)$ is a normalization of 
$\Psi(\mu_bx)$, where $$\mu_b=\Big( 
\frac{1}{2}j_{\frac{N}{2}-1,b}\Big)^2,\quad
\Psi(z)=J_{\frac{N}{2}-1}(2\sqrt{z})(\sqrt{z})^{-\frac{N}{2}+1}.$$
Here $J_{\frac{N}{2}-1}$ is the Bessel function and $j_{\frac{N}{2}-1,b}$
is its $b$-th zero. Since it is known that $\mu_1>0$ and
$\mu_b\sim (b\pi/4)^2$ as $ b\rightarrow +\infty$, we have
$\mu_b\simeq b^2$.

Now we define
$$\tilde{S}(\theta)\tilde{u}=\sum_{a^2\leq\theta, b^2\leq\theta}c_{ab}
\phi_a(t)\psi_b(x)$$
for
$$\tilde{u}=\sum c_{ab}\phi_a(t)\psi_b(x)
\in C_0^{\infty}(]-2T,2T[\times[0,1[),$$
and $S(\theta)u$ be the restriction of
$\tilde{S}(\theta)\tilde{u}$ onto $[0,T]\times[0,1[$.

We claim

\begin{Proposition}
If $0\leq \nu\leq\bar{\nu}$, then
\begin{equation}
\|\tilde{S}(\theta)\tilde{u}\|_{\bar{\nu}}^{\sharp}
\lesssim
\theta^{\bar{\nu}-\nu}\|\tilde{u}\|_{\nu}^{\sharp} \label{G11}
\end{equation}
and
\begin{equation}
\|(I-\tilde{S}(\theta))\tilde{u}\|_{\nu}^{\sharp}
\lesssim \theta^{\nu-\bar{\nu}}
\|\tilde{u}\|_{\bar{\nu}}^{\sharp}
\label{G12}
\end{equation}
\end{Proposition}

Proof.  Note that 
\begin{align*}
(\|\tilde{u}\|_{\nu}^{\sharp})^2&=
\sum_{j+k\leq\nu}
\int_{-2T}^{2T}
\|(-\partial_t^2)^{j}(-\triangle)^{k}\tilde{u}(t,\cdot)\|^2dt \\
&\simeq\sum_{a,b}\sum_{j+k=\nu}
|c_{ab}|^2\lambda_a^{2j}\mu_b^{2k} \\
&\simeq \sum_{a,b}\sum_{j+k=\nu}|c_{ab}|^2(a^jb^k)^4.
\end{align*}
Note that, if $0\leq X, Y\leq 1$ and $0\leq \nu\leq\bar{\nu}$, then
$$\sum_{j+k=\bar{\nu}}X^jY^k\leq
(\bar{\nu}-\nu+1)\sum_{j+k=\nu}X^jY^k. $$
Hence
\begin{align*}
(\|\tilde{S}(\theta)\tilde{u}\|_{\bar{\nu}}^{\sharp})^2&\simeq
\sum_{a^2\leq\theta,b^2\leq\theta}\sum_{j+k=\bar{\nu}}
|c_{ab}|^2(a^jb^k)^4 \\
&=\theta^{2\bar{\nu}}\sum_{a^2\leq\theta,b^2\leq\theta}\sum_{j+k=\bar{\nu}}
|c_{ab}|^2\Big(\frac{a^4}{\theta^2}\Big)^j\Big(\frac{b^4}{\theta^2}\Big)^k \\
&\lesssim\theta^{2\bar{\nu}}
\sum_{a,b}\sum_{j+k=\nu}
|c_{ab}|^2\Big(\frac{a^4}{\theta^2}\Big)^j\Big(\frac{b^4}{\theta^2}\Big)^k \\
&=\theta^{2(\bar{\nu}-\nu)}
\sum_{a,b}\sum_{j+k=\nu}|c_{ab}|^2(a^jb^k)^4 \\
&\simeq \theta^{2(\bar{\nu}-\nu)}(\|\tilde{u}\|_{\nu}^{\sharp})^2,
\end{align*}
that is, (\ref{G11}) holds.

On the other hand, 
$$(\|(I-\tilde{S}(\theta))\tilde{u}\|_{\nu}^{\sharp})^2
\simeq
\sum_{(\theta<a^2)\vee(\theta<b^2)}
\sum_{j+k=\nu}
|c_{ab}|^2(a^jb^k)^4. $$
But 
\begin{align*}
\sum_{\theta<a^2}\sum_{j+k=\nu}|c_{ab}|^2(a^jb^k)^4&=
\theta^{2(\nu-\bar{\nu})}\sum_{\theta<a^2}
\sum_{j+k=\nu}|c_{ab}|^2(\theta^{(\bar{\nu}-\nu)/2}a^jb^k)^4 \\
&\leq\theta^{2(\nu-\bar{\nu})}\sum_{\theta<a^2}
\sum_{j+k=\nu}|c_{ab}|^2(a^{\bar{\nu}-\nu+j}b^k)^4 \\
&\leq\theta^{2(\nu-\bar{\nu})}\sum_{a,b}
\sum_{j+k=\bar{\nu}}|c_{ab}|^2(a^jb^k)^4 \\
&\simeq \theta^{2(\nu-\bar{\nu})}(\|\tilde{u}\|_{\bar{\nu}}^{\sharp})^2.
\end{align*}
We see the same estimate for $\sum_{\theta<b^2}\sum_{j+k=\nu}$.
Hence (\ref{G12}) holds. $\blacksquare$\\

We can define $S_{[0]}(\theta)$ similarly, replacing $N$ by $5$.\\

This conclusions lead us to 

\begin{Proposition}
If $0\leq\nu\leq\bar{\nu}$, then 
\begin{equation}
\|\vec{S}(\theta)\vec{u}\|_{\bar{\nu}}^{(2)}\lesssim \theta^{\bar{\nu}-\nu}
\|\vec{u}\|_{\nu}^{(2)}
\label{G13}
\end{equation}
and
\begin{equation}
\|(I-\vec{S}(\theta))\vec{u}\|_{\nu}^{(2)}
\lesssim
\theta^{\nu-\bar{\nu}}\|\vec{u}\|_{\bar{\nu}}^{(2)}.
\label{G14}
\end{equation}
\end{Proposition}

\section{Estimate of $\mathfrak{F}(\vec{w})$ by $\vec{w}$}

We are going to estimate $\|\mathfrak{F}(\vec{w})\|_{\nu}^{(2)}$ by
$\|\vec{w}\|_{\nu+1}^{(2)}$.\\

Let us consider 
$$\triangle=\triangle_{[1]}=X\frac{d^2}{dX^2}+\frac{N}{2}\frac{d}{dX}, 
\quad D=D_{[1]}=\frac{d}{dX}, \quad \dot{D}=\dot{D}_{[1]}=\sqrt{X}\frac{d}{dX} $$
with $X=1-x$. Let us write $x$ instead of $X=1-x$ and write $\|\cdot\|_k$ instead of $\|\cdot\|_{[1]k}$.\\

Then the first Sobolev's imbedding for functions $u(x)$ reads
\begin{equation}
\|u\|_{L^{\infty}}\lesssim \|u\|_{s_N},
\end{equation}
where
\begin{equation}
 s_N:=\Big[\frac{N}{2}\Big]+1=\min\Big\{ s\in 
\mathbb{N}\  \Big|\  s> \frac{N}{2}\Big\}.
\end{equation}

\begin{Proposition}
We have
\begin{subequations}
\begin{align}
&\|\triangle^mD^ju\|\lesssim \|\triangle^{m+j}u\|, \label{G15a} \\
&\|\dot{D}\triangle^mD^ju\|\lesssim \|\dot{D}\triangle^{m+j}u\|.
\label{G15b}
\end{align}
\end{subequations}
Therefore we can claim
\begin{equation}
\|D^ju\|_k\lesssim \|u\|_{k+2j}. \label{G16}
\end{equation}
\end{Proposition}

Proof. As \cite[Proposition 3]{FE} we start from the formula
$$\triangle^mDu(x)=
x^{-\frac{N}{2}-m-1}\int_0^x
\triangle^{m+1}u(x')(x')^{\frac{N}{2}+m}dx'.$$
Then we can estimate
$$\|\triangle^mDu\|^2=
\int_0^1x^{-\frac{N}{2}-2m-3}\Big(
\int_0^x\triangle^{m+1}u(x')(x')^{\frac{N}{2}+m}dx'
\Big)^2dx$$
by
$$\frac{1}{\frac{N}{2}+m+1}\int_0^1x^{-m-2}
\int_0^x|\triangle^{m+1}u(x')|^2(x')^{\frac{N}{2}+m}dx'dx,$$
using the Schwartz' inequality. It is estimated by
$$\frac{1}{\frac{N}{2}+m+1}\frac{1}{m+1}\|\triangle^{m+1}u\|^2,$$
using the Fubini's change of the order of integrations.
This implies (\ref{G15a}) for $j=1$. Repeating this argument, we get (\ref{G15a}) for $j\geq 2$, too. Let us omit the proof of (\ref{G15b}).$\blacksquare$\\

\begin{Proposition}
We have
\begin{subequations}
\begin{align}
&\|\triangle^m\dot{D}^kD^ju\|\lesssim \|u\|_{2m+k+2j}, \label{G17a}\\
&\|\dot{D}\triangle^m\dot{D}^kD^ju\|\lesssim
\|u\|_{2m+1+k+2j}, \label{G17b}
\end{align}
\end{subequations}
therefore
\begin{equation}
\|\dot{D}^kD^ju\|_n\lesssim\|u\|_{n+k+2j}. \label{G18}
\end{equation}
\end{Proposition}

Proof. First we consider the case when $m=0$. We start from the formula \cite[(B.1)]{FE}:
$$\dot{D}^kDu(x)=
x^{-\frac{N+k}{2}}\int_0^x
\dot{D}^k\triangle u(x')(x')^{\frac{N+k}{2}-1}dx'.$$
Using the Schwartz' inequality and the Fubini's theorem, we can verify
$$\|\dot{D}^kDu\|\leq
\frac{2}{\sqrt{(N+k)k}}\|\dot{D}^k\triangle u\|.$$
On the other hand, using $\dot{D}^2=\triangle-\frac{N-1}{2}D$ and $ [D, \triangle]=D^2$, we can verify that
$$\dot{D}^{2K}\triangle =\sum_{k+\alpha=K+1}C_{k\alpha}\triangle^{\alpha}D^k,$$
with some constants $C_{k\alpha}$. Hence
$$\|\dot{D}^{2K}\triangle u\|\lesssim \|u\|_{2K+2} $$
by Proposition 3, and so on. The conclusion is
$$\|\dot{D}^kD^ju\|\lesssim\|u\|_{k+2j},$$
that is, (\ref{G17a}) holds at least for $m=0$. 

Using $[D,\triangle]=D^2$, we can show
\begin{equation}
D^j\triangle^k=\sum_{\alpha+\beta=k+j}
C_{\alpha\beta}\triangle^{\alpha}D^{\beta}. \label{G19}
\end{equation}
with some constants $C_{\alpha\beta}$.
Then, using (\ref{G19}) and $$[\triangle,\dot{D}]=\frac{N-1}{4}E^{1/2}D,$$
where $E=x^{-1}$, we have
\begin{equation}
[\triangle^m,\dot{D}]=\sum C_{k\alpha\beta
}E^{k+\frac{1}{2}}\triangle^{\alpha}D^{\beta}, \label{G20}
\end{equation}
where the summation runs over
$k+\alpha+\beta=m, \beta\geq 1$. On the other hand, we have
\begin{subequations}
\begin{align}
&s=k\in\mathbb{N} \Rightarrow
\|E^su\|\lesssim
\|D^ku\|, \label{G21a} \\
&s=k+\frac{1}{2}\in \mathbb{N}+\frac{1}{2} \Rightarrow
\|E^su\|\lesssim \|\dot{D}D^ku\|
\label{G21b}.
\end{align}
\end{subequations}
In fact, if $s\in\mathbb{N}/2$, then
\begin{equation}
\|E^su\|\leq\frac{1}{|\frac{N}{4}-s|}\|E^{s-1}Du\|,
\label{G22}
\end{equation}
for the integration by parts leads us to
\begin{align*}
\|E^su\|^2&=\int_0^1u(x)^2x^{-2s+\frac{N}{2}-1}dx \\
&=-\frac{1}{\frac{N}{4}-s}\int_0^1
u\cdot Du\cdot x^{-2s+\frac{N}{2}}dx \\
&\leq\frac{1}{|\frac{N}{4}-s|}\|x^{-s}u\|\|x^{-s+1}Du\|
\end{align*}
provided that $u\in C_0^{\infty}(]0,1[)$. Note that
$N/2 \notin \mathbb{N}$ is supposed by $\neg${\bf (B)} so that
$\frac{N}{4}-s\not=0$ for $s\in\mathbb{N}/2$.
Therefore, if $k+\alpha+\beta=m$, then
\begin{align*}
\|E^{k+\frac{1}{2}}\triangle^{\alpha}D^{\beta}u\|
&\lesssim \|\dot{D}D^k\triangle^{\alpha}D^{\beta}u\| \\
&\lesssim\sum_{\iota+\nu=k+\alpha}\|\dot{D}\triangle^{\iota}D^{\nu}D^{\beta}u\| \\
&\lesssim\sum_{\iota+\gamma=m}\|\dot{D}\triangle^{\iota}D^{\gamma}u\| \\
&\lesssim \|\dot{D}\triangle^mu\|.
\end{align*}
Hence
\begin{equation}
\|\triangle^m\dot{D}u\|\lesssim \|\dot{D}\triangle^mu\|.
\label{G23}
\end{equation}
This can be used to consider the case of odd $k$ in
(\ref{G17a}). If $k=2K$, the identities
$\dot{D}^2=\triangle-\frac{N-1}{2}D$ and (\ref{G19}) are sufficient.
$\blacksquare$\\

\begin{Proposition}
Suppose that
$$\ell_1+\cdots+\ell_p+2(j_1+\cdots+j_p)=n.$$
Then
$$\|(\dot{D}^{\ell_1}D^{j_1}u_1)\cdots
(\dot{D}^{\ell_p}D^{j_p}u_p)\|\lesssim
1+\|u_1\|_n+\cdots \|u_p\|_n, $$
provided that $\|u_{\beta}\|_{2s_N}\lesssim 1\quad \forall \beta$.
\end{Proposition}

Proof. By the Sobolev's imbedding and Proposition 4, (\ref{G18}), we have
\begin{align*}
\clubsuit :=&\|(\dot{D}^{\ell_1}D^{j_1}u_1)\cdots
(\dot{D}^{\ell_p}D^{j_p}u_p)\| \\
&\lesssim
\|u_1\|_{s_N+\ell_1+2j_1}\cdots\|u_{p-1}\|_{s_N+\ell_{p-1}+2j_{p-1}}
\|u_p\|_{\ell_p+2j_p}.
\end{align*}
Suppose $n >s_N$. (Otherwise, we can claim that $\clubsuit \lesssim
\min_{\beta}\|u_{\beta}\|_n$,
provided that $\|u_{\beta}\|_{2s_N}\lesssim 1$.) Suppose that
$\ell_p+2j_p\geq s_N$. (Otherwise, we can assume that
$\ell_{\beta}+2j_{\beta}<s_N$ for $\forall \beta$, and
$\clubsuit \lesssim \min_{\beta}\|u_{\beta}\|_n$, provided that
$\|u_{\beta}\|_{2s_N}\lesssim 1 \forall\beta$. )
by interpolation we have
$$\|u_{\beta}\|_{s_N+\ell_{\beta}+2j_{\beta}}
\lesssim
 \|u_{\beta}\|_{s_N}^{1-\theta_{\beta}}\|u_{\beta}\|_n^{\theta_{\beta}}$$
for $\beta\leq p-1$ with 
$\displaystyle\theta_{\beta}=\frac{\ell_{\beta}+2j_{\beta}}{n-s_N} $,
and
$$\|u_p\|_{\ell_p+2j_p}\lesssim
\|u_p\|_{s_N}^{\frac{n-(\ell_p+2j_p)}{n-s_N}}
\|u_p\|_n^{\frac{\ell_p+2j_p-s_N}{n-s_N}}.
$$
Note that $s_N\leq \ell_p+2j_p\leq n$ and
$$\frac{\ell_p+2j_p-s_N}{n-s_N}=
\frac{n-\sum_{\beta=1}^{p-1}(\ell_{\beta}+2j_{\beta})-s_N}{n-s_N}=
1-\sum_{\beta=1}^{p-1}\theta_{\beta}.$$
Therefore
$$\clubsuit \lesssim
X_1^{\theta_1}\cdots X_{p-1}^{\theta_{p-1}}X_p^{1-(\theta_1+\cdots+\theta_{p-1})}$$
with $X_{\beta}=\|u_{\beta}\|_n$,
provided that
$\|u_{\beta}\|_{s_N}\lesssim 1$. It is easy to show that
$$X_1^{\theta_1}\cdots X_{p-1}^{\theta_{p-1}}X_p^{\theta_p}
\leq X_1+\cdots +X_{p-1}+X_p$$
provided that $0\leq\theta_{\beta}\leq 1, 
\sum_{\beta=1}^p\theta_{\beta}=1$.
$\blacksquare$\\

\begin{Proposition}
Let $F\in C^{\infty}(\mathbb{R}^p)$. Then
$$\|F(u_1,\cdots,u_p)\|_n
\lesssim 1+\|u_1\|_n+\cdots \|u_p\|_n,$$
provided that $\|u_{\beta}\|_{2s_N}\lesssim 1 \quad\forall\beta$.
\end{Proposition}

Proof. As \cite[(B.5)]{FE} we have
$$\triangle^m =\sum_{\mu+j=m}C_{\mu j}\dot{D}^{2\mu}D^j.$$
Therefore
$\triangle^mF(u_1, \cdots, u_p), m\geq 1$ consists of terms of the following types:

(I)  $$(D^{\alpha_1}u_{\beta_1})\cdots(D^{\alpha_q}u_{\beta_q})
\partial_{\beta_q}\cdots\partial_{\beta_1}F $$
with $\alpha_i\geq 1$ for $\forall i$, $\alpha_1+\cdots+\alpha_q=m$. Here
$\partial_{\beta}=\partial/\partial u_{\beta}$.

(II)
$$(\dot{D}^{\gamma_1}D^{\alpha_1}u_{\beta_1})\cdots
(\dot{D}^{\gamma_q}D^{\alpha_q}u_{\beta_q})
\partial_{\beta_q}\cdots\partial_{\beta_1}F$$
with $\alpha_i\geq 1$ for $\forall i$,
$\gamma_i\geq 0$ for $\forall i$,
$\gamma_1+\cdots+\gamma_q=2\mu, \alpha_1+\cdots+\alpha_q=j$,

\noindent$\mu+j=m$.

(III)
$$(\dot{D}^{\gamma_1}D^{\alpha_1}u_{\beta_1})\cdots
(\dot{D}^{\gamma_q}D^{\alpha_q}u_{\beta_q})
(\dot{D}^{\gamma_{q+1}}u_{\delta_1})\cdots
(\dot{D}^{\gamma_{q+r}}u_{\delta_r})
\partial_{\delta_r}\cdots\partial_{\delta_1}
\partial_{\beta_q}\cdots\partial_{\beta_1}F.$$
with $\alpha_i\geq 1$ for $\forall i\leq q$,
$\gamma_k\geq 0$ for $\forall k\leq q$, $\gamma_k\geq 1$ for $\forall k\geq q+1$,

\noindent$\gamma_1+\cdots+\gamma_q+\gamma_{q+1}+\cdots+
\gamma_{q+r}=2\mu$, $\alpha_1+\cdots+\alpha_q=j$, $\mu+j=m$.

(IV)
$$(\dot{D}^{\gamma_1}u_{\delta_1})\cdots
(\dot{D}^{\gamma_r}u_{\delta_r})\partial_{\delta_r}\cdots\partial_{\delta_1}F$$
with
$\gamma_k\geq 1$ for $\forall k\leq r$, $\gamma_1+\cdots+\gamma_r=2m$.

Then we can apply Proposition 5 to each term. $\blacksquare$\\

Now let us consider $\mathfrak{F}(\vec{w})=(\mathfrak{F}_1,
\mathfrak{F}_2)^T$.

Put $\mathfrak{F}_1^{[0]}=\omega(x)\mathfrak{F}_1,
\mathfrak{F}_1^{[1]}=(1-\omega(x))\mathfrak{F}_1$
and so on. 

We are going to estimate
$$(\|\mathfrak{F}(\vec{w})\|_{\nu}^{(2)})^2=
(\|\mathfrak{F}_1\|_{\nu}^{(2)})^2+(\|\mathfrak{F}_2\|_{\nu}^{(2)})^2.$$

Let us observe $\mathfrak{F}_2$. Note that
\begin{align*}
\mathcal{L}y&=\mathcal{L}y^{[0]}+\mathcal{L}y^{[1]} \\
&=-(1-x)\triangle_{[0]}y^{[0]}+\Big(\frac{N}{2}+L_1\Big)Dy^{[0]}+L_0y^{[0]}+ \\
&+x\triangle_{[1]}y^{[1]}+
\Big(-\frac{5}{2}(1-x)+L_1\Big)Dy^{[1]}+L_0y^{[1]}.
\end{align*}
Therefore $\mathfrak{F}_2$ is a smooth function of 
$$1-x, (1-x)^{N/2}, \partial_tv^{[0]}, \partial_tv^{[1]}, 
\partial_ty^{[0]}, \partial_ty^{[1]}, $$
$$y^{[0]}, y^{[1]}, 
Dy^{[0]}, Dy^{[1]}, \triangle_{[0]}y^{[0]}, \triangle_{[1]}y^{[1]}, v^{[0]}, v^{[1]}, Dv^{[0]}, Dv^{[1]}.$$

Let us observe
\begin{align*}
(\|\mathfrak{F}_2\|_{\nu}^{(2)})^2&=
\sum_{j+k\leq\nu}\int_0^T
\|(-\partial_t^2)^j\mathfrak{F}_2(t,\cdot)\|_k^*)^2dt \\
&=\sum_{j+k\leq\nu}\int_0^T(\|(-\partial_t^2)^j\mathfrak{F}_2^{[0]}(t,\cdot)\|_{[0]k}^*)^2dt +\\
&+\sum_{j+k\leq\nu}
\int_0^T(\|(-\partial_t^2)^j\mathfrak{F}_2^{[1]}(t,\cdot)
\|_{[1]k}^*)^2dt.
\end{align*}
When we estimate 
$$\|\triangle_{[1]}^k\mathfrak{F}_2^{[1]}\|_{[1]}=
\|\triangle_{[1]}^k(1-\omega(x))\mathfrak{F}_2\|_{[1]}$$
we see that $\triangle_{[1]}^k(1-\omega)\mathfrak{F}_2$ consists of terms of the form
$$\spadesuit=(\dot{D}^{\gamma_1}D^{\alpha_1}u_{\beta_1})\cdots
(\dot{D}^{\gamma_p}D^{\alpha_p}u_{\beta_p})
\partial_{\beta_p}\cdots\partial_{\beta_1}((1-\omega){F}_2),$$
where $F_2$ is a smooth function of $u_1,\cdots, u_p$. We have to estimate $\|u_{\beta_i}\|_{[1]k}$. Consider the case
in which some of $u_{\beta_i}$ is either $y^{[0]}, Dy^{[0]},
\triangle_{[0]}y^{[0]}, v^{[0]}$ or $
Dv^{[0]}$. E.g., if $u_{\beta_1}=y^{[0]}$, its support is included in
$[0,2/3]$. On the other hand the support of $(1-\omega)F_2$ is included in $[1/3, 1]$. Therefore, if $\chi \in C^{\infty}(\mathbb{R})$ satisfies
$\chi(x)=0$ for $x\leq 1/6$, $0<\chi(x)<1$ for $1/6<x<1/3$, $\chi(x)=1$ for
$1/3\leq x\leq 2/3$, $0<\chi(x)<1$ for $2/3<x<5/6$ and
$\chi(x)=0$ for $5/6\leq x$, we have the equality
$$\spadesuit=(\dot{D}^{\gamma_1}D^{\alpha_1}\chi(x)y^{[0]})
\cdots\partial_{\beta_1}(1-\omega)F_2$$
holds everywhere on $0\leq x\leq 1$. Since the support of $\chi(x)y^{[0]}$ is
included in $[1/6, 5/6]$ , we see 
$$\|\chi(x)y^{[0]}\|_{[1]k}\lesssim
\|\chi(x)y^{[0]}\|_{[0]k}
\lesssim \|y^{[0]}\|_{[0]k}\leq
\|y\|_k.$$
Thus we can claim

\begin{Proposition}
If $\|(1-x)^{N/2}\|_{2s_N+2}<\infty$ and
if $|\|\vec{w}\||_{2s_N+2}\lesssim 1$, then
$$|\|\mathfrak{F}(\vec{w})\||_n\lesssim
1+|\|\vec{w}\||_{n+2},$$
provided that $\|(1-x)^{N/2}\|_n<\infty$.
\end{Proposition}

Note that $\|(1-x)^{N/2}\|_{2s_N+2}<\infty$ if $N>12$.\\

\textbullet \hspace{10mm} $\|\vec{w}\|_{E_1}\lesssim 1$ should imply $|\|\vec{w}\||_{2s_N+2}\lesssim 1$. This requires
that $\|\vec{w}\|_{\mathfrak{b}_E+\mathfrak{r}}^{(2)}\lesssim |\|\vec{w}\||
_{2(\mathfrak{b}_E+\mathfrak{r})} \lesssim 1$
should imply
$|\|\vec{w}\||_{2s_N+2}\lesssim 1$. Hence we require that
\begin{equation}
\mathfrak{b}_E+\mathfrak{r}\geq s_N+2. \label{Sch1}
\end{equation}

\textbullet \hspace{10mm} On the other hand,
$$|\|\mathfrak{F}(\vec{w})\||_n\lesssim
1+|\|\vec{w}|\||_{n+2}$$
implies, by dint of (\ref{X02}), that
$$\|\mathfrak{F}(\vec{w})\|_{\nu}^{(2)}\lesssim 1+\|\vec{w}\|_{\nu+2}^{(2)}$$
with $n=2\nu$.
But it is required that $$\|\mathfrak{F}(\vec{w})\|_{F_j}
\lesssim 1+\|\vec{w}\|_{E_j}.$$
This requires
$$\mathfrak{b}_F+\mathfrak{r}j \leq \mathfrak{b}_E+\mathfrak{r}j -2,$$
or
\begin{equation}
\mathfrak{b}_F\leq \mathfrak{b}_E-2. \label{Sch2}
\end{equation}

\textbullet \hspace{10mm} The estimate should be applied for
$n=2(\mathfrak{b}_F+10\mathfrak{r})$, therefore it should hold that
$\|(1-x)^{N/2}\|_{2(\mathfrak{b}_F+10\mathfrak{r})}<\infty$, that is,
\begin{equation}
\frac{3N}{2}>2(\mathfrak{b}_F+10\mathfrak{r}). \label{Sch6}
\end{equation}

\section{Existence of the inverse 
$\mathfrak{I}(\vec{w})$ of the the Fr\'{e}chet derivative 
$D\mathfrak{F}(\vec{w})$}

We have to analyze the Fr\'{e}chet derivative
$D\mathfrak{F}(\vec{w})$ of the mapping $\mathfrak{F}$
at a given small $\vec{w}=(\tilde{y},\tilde{v})^T$.
For $\vec{h}=(h,k)^T$ we have $D\mathfrak{F}(\vec{w})\vec{h}
=((DF)_1, (DF)_2)^T$, where
\begin{align*}
(DF)_1&=\frac{\partial h}{\partial t}-Jk-\Big((\partial_yJ)v+(\partial_zJ)vx
\frac{\partial}{\partial x}\Big)h, \\
(DF)_2&=\frac{\partial k}{\partial t}+H_1\mathcal{L}h + \\
&+\Big((\partial_yH_1)\mathcal{L}y+\partial_yH_2+
((\partial_zH_1)\mathcal{L}y+\partial_zH_2)
x\frac{\partial}{\partial x}\Big)h + \\
&+\Big((\partial_vH_1)\mathcal{L}y+\partial_vH_2+
(\partial_wH_2)x\frac{\partial}{\partial x}\Big)k.
\end{align*}

Thanks to the assumption {\bf (B3;} N), we see that there are  functions 
$a_{01}, a_{00}, $

\noindent $a_{11}, a_{10}, a_{21}, a_{20}$ of
class $\mathfrak{A}_{(N)}([0,T]\times U^5)$ of
$t, x,y(=y^*+\tilde{y}), Dy, D^2y, v(=v^*+\tilde{v}), Dv $, where $D=\partial/\partial x$, such that
\begin{subequations}
\begin{align}
(DF)_1&=\frac{\partial h}{\partial t}-Jk+
(a_{01}x(1-x)D+a_{00})h, \\
(DF)_2&=\frac{\partial k}{\partial t}+
H_1\mathcal{L}h+
(a_{11}x(1-x)D+a_{10})h+ \nonumber \\
&+(a_{21}x(1-x)D+a_{20})k.
\end{align}
\end{subequations}

Here a smooth function $a(t,x, y_1,\cdots, y_p)$ of 
$[0,T]\times [0,1[\times U\times\cdots\times U$ belongs to
$\mathfrak{A}_{(N)}([0,T]\times U^p)$ if there are analytic functions
$\Phi_0$ on $[0,T]\times]-\delta,\delta[\times U^p$ and
$\Phi_1$ on $[0,T]\times]-\delta,\delta[^2\times U^p$ with $0<\delta \ll 1$ such that
\begin{align*}
&a(t,x, y_1,\cdots )=\Phi_0(t, x, y_1,\cdots)\quad
\mbox{for}\quad 0<x\ll 1,\\
&a(t,x,y_1,\cdots)=\Phi_1(t,1-x, (1-x)^{N/2}, y_1,\cdots)\quad
\mbox{for}\quad 0<1-x\ll 1.
\end{align*}

We put
\begin{align*}
\mathfrak{X}&:=L^2([0,1]; x^{3/2}(1-x)^{N/2-1}dx), \\
\mathfrak{X}^1&:=\{\phi\in\mathfrak{X} |\  
\dot{D}\phi:=\sqrt{x(1-x)}\frac{d\phi}{dx}\in\mathfrak{X}\}, \\
\mathfrak{X}^2&:=\{\phi\in \mathfrak{X}^1 |\  -\Lambda\phi\in\mathfrak{X}\},
\end{align*}
with
\begin{equation}
\Lambda=x(1-x)\frac{d^2}{dx^2}+\Big(\frac{5}{2}(1-x)-
\frac{N}{2}\Big)\frac{d}{dx},
\end{equation}

We claim
\begin{Proposition}
Given $\vec{g}\in C([0,T], \mathfrak{X}^1\times
\mathfrak{X})$, the
equation $D\mathfrak{F}(\vec{w})\vec{h}=\vec{g}$
admits a unique solution $\vec{h}\in C([0,T], \mathfrak{X}^2\times\mathfrak{X}^1)$ such that
$\vec{h}|_{t=0}=\vec{0}$.
\end{Proposition}

Proof. We can rewrite
$$(DF)_2=\frac{\partial k}{\partial t}-H_1\Lambda h
+b_1\check{D}h+b_0h+a_{21}\check{D}k+a_{20}k,
$$
where
\begin{align*}
\check{D}&:=x(1-x)D, \\
b_1&:=H_1\frac{L_1}{x(1-x)}+a_{11}, \\
b_0&:=H_1L_0+a_{10}.
\end{align*}
Of course $b_1, b_0$ are also analytic on a neighborhood of
$[0,1]\times\{0\}\times\cdots$.
Then we can write the equation
$D\mathfrak{F}(\vec{w})\vec{h}=\vec{g}(=(g_1,g_2)^T)$ as
\begin{equation}
\frac{\partial}{\partial t}
\begin{bmatrix}
h \\
k
\end{bmatrix}
+
\begin{bmatrix}
\mathfrak{a}_1 & -J \\
\mathcal{A} & \mathfrak{a}_2
\end{bmatrix}
\begin{bmatrix}
h \\
k
\end{bmatrix}
=
\begin{bmatrix}
g_1 \\
g_2
\end{bmatrix}.
\label{Eq16}
\end{equation}
Here
\begin{align*}
\mathfrak{a}_1&=a_{01}\check{D}+a_{00}, \\
\mathfrak{a}_2&:=a_{21}\check{D}+a_{20}, \\
\mathcal{A}&:=-H_1\Lambda + b_1\check{D}+b_0.
\end{align*}
The standard calculation gives
\begin{align}
&\frac{1}{2}\frac{d}{dt}\Big(\|k\|^2+((H_1/J)\dot{D}h|\dot{D}h)\Big) +\nonumber\\
&+(\beta_1\dot{D}h|\dot{D}h)+
(\beta_2\dot{D}h|h)+(\beta_3\dot{D}h|k)+
+(\beta_4h|k)+(\beta_5k|k)= \nonumber\\
&=((H_1/J)\dot{D}h|\dot{D}g_1)+(k|g_2), \label{Eneqy}
\end{align}
where $\dot{D}=\sqrt{x(1-x)}D$ and
\begin{align*}
\beta_1&=-\frac{1}{4}(3+(N+3)x+2\check{D})(H_1/J)a_{01}-\frac{1}{2}
\frac{\partial(H_1/J)}{\partial t} +
(H_1/J)(\check{D}a_{01}+a_{00}), \\
\beta_2&=(H_1/J)\dot{D}a_{00}, \\
\beta_3&=-(H_1/J)\dot{D}J+\dot{D}H_1+
\sqrt{x(1-x)}(b_1+a_{21}), \\
\beta_4&=b_0,\qquad \beta_5=a_{20}.
\end{align*}

We assume that $\beta_j \in C([0,T]\times [0,1])$ for $j=1,2,\cdots, 5$.\\

Of course $(\cdot|\cdot)$ and $\|\cdot\|$ stand for the inner product and the norm of the Hilbert space $\mathfrak{X}$. We have used the following
formulas:\\

{\bf Formula 1}: {\it If $\phi\in\mathfrak{X}^2, \psi\in\mathfrak{X}^1,
\alpha\in C^1([0,1])$, then
\begin{equation}
(-\alpha\Lambda\phi|\psi)=
(\alpha\dot{D}\phi|\dot{D}\psi)+
((D\alpha)\check{D}\phi|\psi). \label{F1}
\end{equation}}\\

{\bf Formula 2}: {\it If $\phi \in\mathfrak{X}^2$ and $\alpha\in
C^1([0,1])$, then
\begin{equation}
(\alpha\dot{D}\phi|\dot{D}\check{D}\phi)=
(\alpha^*\dot{D}\phi|\dot{D}\phi), 
\label{F2}
\end{equation}
with
$$\alpha^*=-\frac{1}{4}(3+(N+3)x+2\check{D})\alpha.$$ }\\

Since $\vec{w} =(y^*+\tilde{y}, v^*+\tilde{v})^T$ is
supposed to be small,
we can assume
$$\frac{1}{M_0}<J<M_0,\qquad \frac{1}{M_0}<H_1<M_0 $$
with a constant $M_0$ independent of  $\vec{w}$ thanks to the assumption {\bf (B2)}. Now the energy
$$\mathcal{E}:=\|k\|^2+((H_1/J)\dot{D}h|\dot{D}h) $$
enjoys the inequality
$$\frac{1}{2}\frac{d\mathcal{E}}{dt}\leq
M(\|\vec{h}\|_{\mathfrak{H}}^2+
\|h\|_{\mathfrak{H}}\|\vec{g}\|_{\mathfrak{H}}),
$$
where $\mathfrak{H}=\mathfrak{X}^1\times\mathfrak{X}$ and
$$
\|(\phi,\psi)^T\|_{\mathfrak{H}}^2=\|\phi\|_{\mathfrak{X}^1}^2+
\|\psi\|_{\mathfrak{X}}^2
=\|\phi\|^2+\|\dot{D}\phi\|^2+\|\psi\|^2$$
and
$$M=\sum_{j=1}^5\|\beta_j\|_{L^{\infty}}+(M_0)^2+1.
$$
Since $\mathcal{E}$ is equivalent to $\|k\|^2+\|\dot{D}h\|^2$, the Gronwall
argument and application of the Kato's theory (\cite{Kato11})
deduce the conclusion. Here $\|h\|$ should be estimated by $\mathcal{E}$
as follows: The first component of (\ref{Eq16}) implies
$$h(t)=\int_0^t(-a_{01}\check{D}h-a_{00}h+Jk+g_1)(t')dt', $$
therefore
$$\|h(t)\|\leq C\int_0^t\|h(t')\|dt'+
\int_0^t(C\mathcal{E}(t')+\|g_1(t')\|)dt',$$
where
$C=\max(\|a_{00}\|_{L^{\infty}},\|a_{01}\|_{L^{\infty}}M_0^2+ M_0)$, which implies, through the Gronwall's argument,
$$\|h(t)\|\leq\int_0^t
(e^{C(t-t')}-1)(C\mathcal{E}(t')+\|g_1(t')\|)dt'.$$
As the result the solution enjoys
$$\|\vec{h}(t)\|_{\mathfrak{H}}
\leq C\int_0^t\|\vec{g}(t')\|_{\mathfrak{H}}dt'.
$$
$\blacksquare$\\

Here, in order to make sure, let us sketch proofs of {\bf Formula 1},
{\bf Formula 2}, and  (\ref{Eneqy}).\\

Proof of (\ref{F1}): if $\psi\in \mathfrak{X}^1$, then 
$$\psi(1/2)+\int_{1/2}^x\frac{\dot{D}\psi(x')}{\sqrt{x'(1-x')}}dx'$$
implies
$$|\psi(x)|\leq C x^{-3/4}(1-x)^{-N/4+1/2},$$
and, if $\phi\in\mathfrak{X}^2$, then
\begin{align*}
x^{5/2}(1-x)^{N/2}\frac{d\phi}{dx}&=
x^{5/2}(1-x)^{N/2}\frac{d\phi}{dx}\Big|_{x=1/2} +\\
&-\int_{1/2}^x
\Lambda\phi(x')x'^{3/2}(1-x')^{N/2-1}dx'
\end{align*}
implies
$$
\Big|\frac{d\phi}{dx}\Big|\leq C
x^{-5/4}(1-x)^{-N/4}.$$
Actually the finite constant
$$x^{5/2}(1-x)^{N/2}D\phi|_{x=1/2}+
\int_0^{1/2}\Lambda\phi(x')x'^{3/2}(1-x')^{N/2-1}dx'$$
should vanish in order to $\dot{D}\phi \in\mathfrak{X}$ and so on. 
Therefore the boundary terms in the integration by parts vanish
at $x=+0, 1-0$ and we get the desired equality.\\

Proof of (\ref{F2}): We see
\begin{align*}
(\alpha\dot{D}\phi|\dot{D}\check{D}\phi)&=
\int_0^1\alpha x(1-x)(D\phi)D(x(1-x)D\phi)x^{3/2}(1-x)^{N/2-1}dx \\
&=I +(\alpha(1-2x)\dot{D}\phi|\dot{D}\phi),
\end{align*}
where
\begin{align*}
I&:=\int_0^1\alpha(x(1-x))^2(D\phi)(D^2\phi)x^{3/2}(1-x){N/2-1}dx \\
&=\int_0^1\frac{\alpha}{2}D(D\phi)^2x^{7/2}(1-x)^{N/2+1}dx \\
&=-\int_0^1
D\Big(\frac{\alpha}{2}x^{7/2}(1-x)^{N/2+1}\Big)(D\phi)^2dx.
\end{align*}
Here the integration by parts has been done by using
$$|D\phi|\leq C x^{-5/4}(1-x)^{-N/4}$$
which holds for $\phi\in\mathfrak{X}^2$. 
Then we see
$$I=-\Big(\Big(\frac{\check{D}\alpha}{2}+
\frac{\alpha}{2}\Big(\frac{7}{2}(1-x)+\Big(\frac{N}{2}+1\Big)x\Big)\Big)
\dot{D}\phi\Big|\dot{D}\phi\Big),$$
and get (\ref{F2}). \\

Proof of (\ref{Eneqy}): multiplying the second component of the equation
(\ref{Eq16}) by $k$ and integrating it, we get
\begin{align*}
\frac{1}{2}\frac{d}{dt}
\|k\|^2-(H_1\Lambda h|k)&+(b_1\check{D}h|k)+(b_0h|k) + \\
&+(a_{21}\check{D}h|k)+(a_{20}k|k)=(g_2|k).
\end{align*}
By {\bf Formula 1} we see
\begin{align}
\frac{1}{2}\frac{d}{dt}
\|k\|^2
+(H_1\dot{D}h|\dot{D}k)&
+((DH_1)\check{D}h|k)
+(b_1\check{D}h|k)+(b_0h|k) + \nonumber\\
&+(a_{21}\check{D}h|k)+(a_{20}k|k)=(g_2|k). \label{*}
\end{align}
On the other hand, operating $\dot{D}$ on the first component of
(\ref{Eq16}), we get
\begin{align*}
\dot{D}k&=\frac{1}{J}\partial_t(\dot{D}h)+
\frac{a_{01}}{J}\dot{D}\check{D}h+ \\
&+\frac{1}{J}(\check{D}a_{01}+a_{00})\dot{D}h+
\frac{\dot{D}a_{00}}{J}h-\frac{\dot{D}J}{J}h-\frac{\dot{D}J}{J}k-
\frac{1}{J}\dot{D}g_1.
\end{align*}
Inserting this into the second term of the left-hand
side of (\ref{*}), we get
\begin{align*}
&\frac{1}{2}\frac{d}{dt}\|k\|^2+\frac{1}{2}\frac{d}{dt}
((H_1/J)\dot{D}h|\dot{D}h)+
((H_1a_{01}/J)\dot{D}h|\dot{D}\check{D}h)+ \\
&+(((1/2)(H_1/J)_t+(H_1/J)(\check{D}a_{01}+a_{00}))\dot{D}h|\dot{D}h)+
((H_1/J)(\dot{D}a_{00})\dot{D}h|h) + \\
&-((H_1\dot{D}J/J)\dot{D}h|k)-
((H_1/J)\dot{D}h|\dot{D}g_1) + \\
&+((b_1+a_{21})\check{D}h|k)+(b_0h|k)+(a_{20}k|k)=(g_2|k).
\end{align*}
Applying {\bf Formula 2} to the third term of the left-hand
side, we get the desired (\ref{Eneqy}).\\

Of course, when we apply the Nash-Moser(-Schwartz) theorem, we take 
$\mathfrak{I}(\vec{w}): \vec{g} \mapsto \vec{h}$, where $\vec{h}$ is the solution of 
$$D\mathfrak{F}(\vec{w})\vec{h}=\vec{g}, \qquad \vec{h}|_{t=0}=\vec{0}.
$$

\section{Tame estimate of the inverse 
$\mathfrak{I}(\vec{w})$ of the the Fr\'{e}chet derivative 
$D\mathfrak{F}(\vec{w})$}

Now we investigate the equation $D\mathfrak{F}(\vec{w})\vec{h}=
\vec{g}$ on $0\leq t\leq T$. The equation can be written as:
\begin{equation}
\frac{\partial}{\partial t}
\begin{bmatrix}
h^{[\mu]} \\
k^{[\mu]}
\end{bmatrix}+
\begin{bmatrix}
\mathfrak{a}_1^{[\mu]} & -J \\
\mathcal{A}^{[\mu]} & \mathfrak{a}_2^{[\mu]}
\end{bmatrix}
\begin{bmatrix}
h^{[\mu]} \\
k^{[\mu]}
\end{bmatrix}
=\begin{bmatrix}
g_1^{\mu]} \\
g_2^{\mu]}
\end{bmatrix}
+(-1)^{\mu}
\begin{bmatrix}
c_{11} & 0 \\
\mathfrak{c}_{21} & c_{22}
\end{bmatrix}
\begin{bmatrix}
h^{[1-\mu]} \\
k^{[1-\mu]}
\end{bmatrix},
\quad \mu=0,1,
\label{Eq301}
\end{equation}
where
\begin{align*}
&\mathfrak{a}_1^{[\mu]}=a_{01}\check{D}+a_{00}^{[\mu]}=a_{01}\check{D}+
a_{00}-(-1)^{\mu}a_{01}\check{D}\omega, \\
&\mathfrak{a}_2^{[\mu]}=a_{21}\check{D}+a_{20}^{[\mu]}=
a_{21}\check{D}+a_{20}-(-1)^{\mu}a_{21}\check{D}\omega, \\
\mathcal{A}^{[\mu]}&=-H_1\Lambda+
(b_1+(-1)^{\mu}2H_1(D\omega))\check{D}+b_0+
(-1)^{\mu}(H_1\Lambda-b_1\check{D})\omega, \\
c_{11}&=a_{01}\check{D}\omega, \\
\mathfrak{c}_{21}&=-2H_1(D\omega)\check{D}+
b_1\check{D}\omega-H_1\Lambda\omega, \\
c_{22}&=a_{21}\check{D}\omega.
\end{align*}

We can write
$$\mathcal{A}^{[\mu]}=-b_2^{[\mu]}\triangle_{[\mu]}+b_1^{[\mu]}\check{D}_{[\mu]}+b_0^{[\mu]},
$$
where
\begin{align*}
&\check{D}_{[0]}=x\frac{\partial}{\partial x},
\qquad
\check{D}_{[1]}=X\frac{\partial}{\partial X}\quad\mbox{with}\quad X=1-x, \\
& b_2^{[0]}=H_1\cdot (1-x), \qquad b_1^{[0]}=\frac{N}{2}H_1+
(b_1+2H_1(D\omega))(1-x), \\
&b_2^{[1]}=H_1\cdot x, \qquad b_1^{[1]}=\frac{5}{2}H_1-
(b_1-2H_1(D\omega))x, \\
&b_0^{[\mu]}=b_0-(-1)^{\mu}(H_1\Lambda -b_1\check{D})\omega.
\end{align*}
First we consider the component $\vec{h}=\vec{h}^{[\mu]}, \mu=0,1$, which satisfies
$$\frac{\partial\vec{h}}{\partial t}+
\mathfrak{A}\vec{h}=\vec{f},$$
with
\begin{align*}
&\mathfrak{A}=\mathfrak{A}_{[\mu]}
=\begin{bmatrix}
 \mathfrak{a}_1^{[\mu]} & -J \\
\mathcal{A}^{[\mu]} & \mathfrak{a}_2^{[\mu]}
\end{bmatrix}, \\
&\vec{f}=\vec{f}^{[\mu]}=
\begin{bmatrix}
f_1^{[\mu]} \\
f_2^{[\mu]}
\end{bmatrix}
= 
\begin{bmatrix}
g_1^{[\mu]} \\
g_2^{[\mu]}
\end{bmatrix}
+(-1)^{\mu}
\begin{bmatrix}
c_{11} & 0 \\
\mathfrak{c}_{21} & c_{22}
\end{bmatrix}
\begin{bmatrix}
h^{[1-\mu]} \\
k^{[1-\mu]}
\end{bmatrix}.
\end{align*}

{\bf We consider $\mu=1$. } Let us write $x$ instead of $X=1-x$.\\

First we prepare the second Sobolev's imbedding:

\begin{Proposition}
If $s\in \mathbb{N}, s<\displaystyle \frac{N}{2}, 
\frac{1}{2}-\frac{s}{N}\leq \frac{1}{p}\leq \frac{1}{2}$,
then 
$\|u; L^p\|\lesssim \|u\|_s$. 

(Note that
$s<N/2 \Leftrightarrow s<s_N$ for $N/2\not\in \mathbb{N}$.)
\end{Proposition}

Here we denote the $L^p$-norm with respect to the measure
$x^{N/2-1}dx$ by
$$\|u; L^p\|:=\Big(\int_0^1
|u(x)|^px^{\frac{N}{2}-1}dx\Big)^{1/p}.$$

Proof. We may assume that $u\in C_0^{\infty}([0,1[)$, by dint of the extention technique. Supposet that $s\geq 1, 
\displaystyle \frac{1}{2}-\frac{s}{N}<\frac{1}{p}\leq \frac{1}{2}$. Let
$$u(x)=\sum_{n=1}^{\infty}c_n\psi_n(x),$$
where
\begin{align*}
&\psi_n(x)=\frac{\Psi_{\nu}(\lambda_nx)}{\|\Psi_{\nu}(\lambda_nx)\|},\\
&\Psi_{\nu}\Big(\frac{r^2}{4}\Big)=J_{\nu}(r)\Big(\frac{r}{2}\Big)^{-\nu},\qquad
\nu=\frac{N}{2}-1, \\
&\lambda_n=\Big(\frac{j_{\nu,n}}{2}\Big)^2,
\end{align*}
Here $J_{\nu}$ is the Bessel function, and $j_{\nu,n}$ is its
$n$-th positive zero.. See \cite[Appendix A]{FE}. We know
$$\|\Psi_{\nu}(\lambda_nx)\|^{-1}\simeq n^{\frac{N-1}{2}}.$$
Since 
$$\int_0^1|\Psi_{\nu}(\lambda_nx)|^px^{\nu}dx=\frac{1}{\lambda_n^{\nu+1}}I,$$
with
$$I=2^{\nu(p-2)-1}\int_0^{j_{\nu, n}}
|J_{\nu}(r)|^pr^{-\nu(p-2)+1}dr,$$
and since
$$J_{\nu}(r)
=
\sqrt{\frac{2}{\pi r}}
\cos\Big(r-\frac{\nu}{2}\pi-\frac{\pi}{4}\Big)+
O\Big(\frac{1}{r^{3/2}}\Big), $$
we see

(1): $I\lesssim 1$, provided that 
$$-\frac{N-1}{2}(p-2)<-1;$$

(2): $I\lesssim \log n$, provided that
$$-\frac{N-1}{2}(p-2)=-1;$$

(3): $I\lesssim n^{-\frac{N-1}{2}(p-2)+1}$, provided that
$$-\frac{N-1}{2}(p-2)>-1.$$

Since $\lambda_n\simeq n^2$, we have
$$\|\psi_n; L^p\|\lesssim n^{\frac{N-1}{2}-\frac{N}{p}}$$
for the case (1); 
$$\|\psi_n;L^p\|\lesssim n^{\frac{N-1}{2}+(-N+\epsilon)\frac{1}{p}},
\quad\mbox{with}\quad 0<\epsilon \ll 1$$
for the case (2); and
$$\|\psi_n; L^p\|\lesssim 1$$
for the case (3). Therefore, keeping in mind that
$\|u\|_s\simeq \sqrt{|c_n|^2\lambda_n^s}$,
we have
$$\|u; L^p\|\lesssim\|u\|_s
\sqrt{\sum n^{2(\frac{N-1}{2}-\frac{N}{p}-s)}}$$
for the case (1); 
$$\|u; L^p\|\lesssim \|u\|_s
\sqrt{\sum n^{2(\frac{N-1}{2}-(N-\epsilon)\frac{1}{p}-s)}}$$
for the case (2); and
$$\|u;L^p\|\lesssim\|u\|_s
\sqrt{n^{-2s}}$$
for the case (3).
Then we see
$$2\Big(\frac{N-1}{2}-\frac{N}{p}-s\Big)
<2\Big(\frac{N-1}{2}-(N-\epsilon)\frac{1}{p}-s\Big)<-1$$
 and
$ -2s<-1$, provided that
$s\geq 1$ and $\frac{1}{2}-\frac{s}{N}<\frac{1}{p}, \epsilon\ll 1$. So,
in every case $\sqrt{\sum\cdots}<\infty$. $\blacksquare$\\

This implies:

\begin{Proposition}
If $s_1, s_2, k\in \mathbb{N}$ satisfy
$s_1\geq k, s_2\geq k, s_1+s_2\geq s_N+k$, then it hols that
$$\|f\cdot g\|_k\lesssim \|f\|_{s_1}\|g\|_{s_2}.$$

(Note that $s_1+s_2>\frac{N}{2}+k\Leftrightarrow
s_1+s_2\geq s_N+k$ for $\frac{N}{2}\not\in\mathbb{N}$.)
\end{Proposition}

Proof. First consider the case with $k=0$. Suppose
$s_1<N/2, s_2<N/2$. H\"{o}lder's inequality
gives
$$\|f\cdot g\|\leq \|f;L^p\|\|g; L^q\|, $$
provided that $\displaystyle \frac{1}{p}+\frac{1}{q}=\frac{1}{2}$. Choose $p$ such that
$$\frac{1}{2}-\frac{s_1}{N}\leq \frac{1}{p}\leq \frac{s_2}{N}\quad \Big(<
\frac{1}{2}\Big). $$
This possible, since $s_1+s_2>N/2$. Then Proposition 8 gives
$\|f; L^p\|\lesssim \|f\|_{s_1}$, and
$$\frac{1}{2}-\frac{s_2}{N}\leq \frac{1}{q}=\frac{1}{2}-\frac{1}{p}
\leq \frac{s_1}{N},$$
which gives $\|g;L^q\|\lesssim \|g\|_{s_2}$. That is done. If $s_1\geq N/2$, 
then $s_1>N/2$, and 
$$\|f;L^{\infty}\|\leq \|f\|_{s_1},\quad
\|g: L^2\|\lesssim \|g\|_{s_2}, \quad 
\|f\cdot g\|\lesssim
\|f;L^{\infty}\|\|g;L^2\|$$ will work for $s_2\geq 0$. Thus the proof for $k=0$ is done.

Next consider the case $k=2m\geq 2$ with $m\in\mathbb{N}$. We have to estimate
$\|\triangle^m(f\cdot g)\|$. But, according to
\cite[(8.5)]{FE}, $\triangle^m(f\cdot g)$ is a linear conbination of terms
$$\spadesuit =(\dot{D}^{\ell'}D^{j'}f)\cdot
(\dot{D}^{\ell}D^jg)$$
with $\ell'+\ell+2(j'+j)=2m$. In order to estimate $\|\spadesuit\|$, we put
$t_1=s_1-(\ell'+2j'), t_2=s_2-
(\ell +2j)$.
Then $s_1\geq 2m, s_2\geq 2m, s_1+s_2\geq s_N+2m$ imply
$t_1\geq 0, t_2\geq 0, t_2+t_2\geq s_N$. Therefore the inequality for $k=0$ can be applied to get
$$\|\spadesuit\|\lesssim
\|\dot{D}^{\ell'}D^{j'}f\|_{t_1}
\|\dot{D}^{\ell}D^jg\|_{t_2}
\lesssim \|f\|_{s_1}\|g\|_{s_2}.$$
We omit the rest of the proof. $\blacksquare$.\\

Let us deduce an elliptic estimate. Let $n,\sigma \in\mathbb{N}$
satisfy $s_N+2\leq\sigma, n+2\leq \sigma$.\\

We know
\begin{equation}
\|[\triangle^m, \mathcal{A}]u\|\lesssim C(m,u),
\end{equation}
with
\begin{align}
C(m,u):=&\sum(
\|\triangle^{\gamma}D^{\delta}b_i)(\dot{D}\triangle^ju)\|+
\|(\triangle^{\gamma}D^{\delta}b_i)(\triangle^ju)\|+\nonumber \\
&
\|(\check{D}\triangle^{\gamma}D^{\delta}b_i)(\triangle^ju)\|).
\end{align}
Here the summation runs over
$$i=0,1,2, \qquad \gamma+\delta\leq k+1,\qquad k+j=m.$$
See \cite[Appendix D]{FE}. Similarly
\begin{equation}
\|\dot{D}[\triangle^m, \mathcal{A}]u\|\lesssim C^{\sharp}(m,u),
\end{equation}
with
\begin{equation}
C^{\sharp}(m,u):=C(m,u)+
\sum \|(\triangle^{\gamma}D^{\delta}b_i)(\triangle^{j+1}u)\|+
\|\dot{D}\triangle^{\gamma}D^{\delta}b_i)(\triangle^ju)\|.
\end{equation}
The range of the summation is the same as above.

Using $\triangle^m[\triangle, \mathcal{A}]=
[\triangle^{m+1},\mathcal{A}]-
[\triangle^m, \mathcal{A}]\triangle$ and
noting that $C(m, \triangle u)\leq C(m+1, u)$, we have
$$\|\triangle^m[\triangle,\mathcal{A}]u\|\lesssim C(m+1,u).$$
We can estimate $\|\dot{D}\triangle^m[\triangle, \mathcal{A}]u\|$
similarly, and we see
\begin{equation}
\|[\triangle,\mathcal{A}]u\|_n\lesssim K(n+2, u),
\end{equation}
where
\begin{equation}
K(n,u):=
\begin{cases} C(m,u)\quad\mbox{for}\quad n=2m \\
C^{\sharp}(m,u)\quad\mbox{for}\quad
n=2m+1
\end{cases}
\end{equation}

Let us estimate $C(m,u)$. For example, consider
$$\clubsuit=\|(\triangle^{\gamma}D^{\delta}b_i)(\dot{D}\triangle^ju)\|,$$
where $i=0,1,2, \gamma+\delta\leq k+1, k+j=m$. Suppose that
$s_N+2\leq\sigma, n+2=2m+2\leq\sigma$. Put
$s_1=\sigma-2-2k, s_2=2k$. Then
$s_1\geq 0, s_2\geq 0, s_1+s_2=\sigma-2\geq s_N$.
Therefore Proposition 10 can be applied to get
\begin{align*}
\clubsuit&\lesssim \|\triangle^{\gamma}D^{\delta}b_i\|_{s_1}
\|\dot{D}\triangle^ju\|_{s_2} \\
&\lesssim \|\vec{b}\|_{s_1+2k+2}\|u|_{s_2+2j+1}=\|\vec{b}\|_{\sigma}
\|u\|_{2m+1}.
\end{align*}
Here
$$\vec{b}=\vec{b}^{[\mu]}=(b_0^{[\mu]}, b_1^{[\mu]}, b_2^{[\mu]}).$$

In this way we can verify
\begin{equation}
K(n,u)\lesssim \|\vec{b}\|_{\sigma}\|u\|_{n+1}, \label{H00}
\end{equation}
provided that $s_N+2\leq\sigma, n+2\leq\sigma$.

Then we have

\begin{Proposition}
If $n, \sigma\in\mathbb{N}$ satisfy
$s_N+2\leq\sigma, n+2\leq \sigma$, then it holds that
\begin{equation}
\|u\|_{n+2}\lesssim \|\mathcal{A}u\|_n+\|u\|_1, \label{H05}
\end{equation}
provided that $\|\vec{b}\|_{\sigma}\lesssim 1$.
\end{Proposition}

Proof. $\triangle u=-\frac{1}{b_2}(\mathcal{A}u-b_1\check{D}u-b_0u)$ implies
\begin{align*}
&\|u\|_2\lesssim \|\mathcal{A}u\|+\|u\|_1,\\
&\|u\|_3\lesssim \|\mathcal{A}u\|_1+\|u\|_1.
\end{align*}
See \cite[pp.83-84]{FE}. Keeping in mind that
$\|u\|_1\lesssim K(n,u)$ for $n\geq 0$, we can verify
$$\|u\|_{n+2}\lesssim \|\mathcal{A}u\|_n+K(n,u)$$
by induction, using the fact
$K(n-2,\triangle u)\leq K(n, u)$ and the interpolation
$$\|u\|_{n+2}\lesssim \|u\|_n+\|\triangle u\|_n.$$
Then (\ref{H00}) implies
$$\|u\|_{n+2}\lesssim \|\mathcal{A}u\|_n+\|u\|_{n+1},
$$
provided that $\|\vec{b}\|_{\sigma}\lesssim 1$. By induction, we can 
replace $\|u\|_{n+1}$ in the left-hand side by $\|u\|_1$.
$\blacksquare$.\\

This implies

\begin{Proposition}
If $n,\sigma\in\mathbb{N}$ satisfy $s_N+2\leq\sigma,
n+2\leq\sigma$, then  it holds that
\begin{equation}
\|\vec{u}\|_{n+1}\lesssim
\|\mathfrak{A}\vec{u}\|_n+\|\vec{u}\|_1,
\end{equation}
provided that $\|\vec{a}\|_{\sigma}\lesssim 1$.
\end{Proposition}
Here
$$\vec{a}=\vec{a}^{[\mu]}=(b_0^{[\mu]}, b_1^{[\mu]}, b_2^{[\mu]}, a_{01}, a_{00}, a_{21}, a_{20}, J).$$

In order to derive Proposition 12 from Proposition 11, it is sufficient to
note that
\begin{equation}
\|\phi\cdot u\|_n\lesssim \|\phi\|_{\sigma'}\|u\|_n,
\end{equation}
provided that $s_N\leq\sigma', n\leq \sigma'$.
We omit the proof. \\

Inversely

\begin{Proposition}
If $n,\sigma\in\mathbb{N}$ satisfy $s_N+2\leq\sigma,
n+2\leq\sigma$, then it holds that
\begin{equation}
\|\mathfrak{A}\vec{u}\|_n\lesssim \|\vec{u}\|_{n+1}
\end{equation}
for $0\leq t\leq \tau$, provided that
$|\vec{a}; \tau, \sigma\|\lesssim 1$.
\end{Proposition}

Here 
we recall the definition
$$|\phi; \tau,\nu\|=|\phi; \tau,\nu\|_{[\mu]}:=
\sup_{0\leq t\leq \tau}\sum_{\iota+\kappa\leq\nu}
\|\partial_t^{\iota}\phi(t,\cdot)\|_{[\mu]\kappa}.
$$
The proof can be reduced to
$$\|\mathcal{A}u\|_n\lesssim\|u\|_{n+2}+K(n,u).$$
Let us omit the details.\\

Let $n,\sigma\in\mathbb{N}$ satisfy
$s_N+1\leq\sigma, n\leq \sigma$. Then we claim
\begin{equation}
\|[\partial_t^j,\mathcal{A}]u\|_k
\lesssim
|\vec{b}; \tau, \sigma\|\cdot
|u;\tau, n+1\|
\end{equation}
for $j+k=n, 0\leq t\leq \tau$.

In fact, let us consider, e.g., 
$$[\partial_t^j,b_2\triangle]u=\sum C_{\alpha\beta}
(\partial_t^{\alpha}b_2)(\partial_t^{\beta}\triangle u).$$
The summation runs over
$$\alpha+\beta=j, \qquad \alpha\geq 1.$$
Put $t_1=\sigma-\alpha, t_2=n-1-\beta=\alpha-1$. Then
$t_1-k=\sigma-n+\beta\geq\beta\geq 0, t_2-k=\alpha-1\geq 0,
t_1+t_2-s_N-k=\sigma-1-s_N\geq 0$. Proposition 10 can be applied to get
\begin{align*}
\|(\partial_t^{\alpha}b_2)(\partial_t^{\beta}\triangle u\|_k&
\lesssim\|\partial_t^{\alpha}b_2\|_{\sigma-\alpha}
\|\partial_t^{\beta}\triangle u\|_{n-1-\beta} \\
&\leq|b_2;\tau, \sigma\|\cdot|u;\tau,n+1\|,
\end{align*}
and so on. This observation implies

\begin{Proposition}
If $s_N+2\leq\sigma, n+1\leq\sigma$, then it holds that
\begin{equation}
\|[\partial_t^j,\mathfrak{A}]\vec{u}\|_k
\lesssim |\vec{a};\tau, \sigma\|\cdot
|\vec{u}; \tau, n\|
\end{equation}
for $j+k=n, 0\leq t\leq\tau$.
\end{Proposition}

Now we are going to find estimates of the solution 
$\vec{h}=\vec{h}^{[\mu]}$ of the problem
$$\frac{\partial \vec{h}}{\partial t}+\mathfrak{A}\vec{h}=
\vec{f}, \qquad
\vec{h}|_{t=0}=\vec{0},
$$
where $\mathfrak{A}=\mathfrak{A}_{[\mu]}, \vec{f}=
\vec{f}^{[\mu]}$.\\

Consider $\mu=1$. Recall that a solution of
$$\frac{\partial\vec{H}}{\partial t}+\mathfrak{A}\vec{H}=\vec{F}$$
enjoys the energy estimate
$$\|\vec{H}(t)\|\lesssim
\|\vec{H}(0)\|+\int_0^t\|\vec{F}(t')\|dt'.$$ \\

We put
\begin{align*}
&X(\vec{H};j,k):=\|\partial_t^j\vec{H}\|_k, \qquad
X_0(\vec{H};j,k):=X(\vec{H};j,k)|_{t=0}, \\
&Z(\vec{H};n):=\sum_{j+k=n}X(\vec{H};j,k), \qquad Z_0(\vec{H};n):=Z(\vec{H};n)|_{t=0}, \\
&W(\vec{H}; n):=\sum_{j+k\leq n}X(\vec{H};j,k)=\sum_{\nu\leq n}
Z(\vec{H};\nu), \qquad W_0(\vec{H};n):=W(\vec{H};n)|_{t=0}.
\end{align*}

First we note
\begin{equation}
Z(\vec{h};n+1)=Z(\partial_t\vec{h};n)+\|\vec{h}\|_{n+1}.
\end{equation}
Suppose that $s_N+2\leq\sigma, n+2\leq\sigma, |\vec{a};\tau,\sigma\|\lesssim 1, 0\leq t\leq\tau$. We claim
\begin{equation}
Z(\vec{h}; n+1)\lesssim
\|\partial_t^{n+1}\vec{h}\|+
\sum_{0\leq \nu\leq n}\|\partial_t^{\nu}\vec{h}\|+\|f\|_n.
\label{H11}
\end{equation}
(\ref{H11}) is true for $n=0$, since
\begin{align*}
Z(\vec{h};1)&=\|\partial_t\vec{h}\|+\|\vec{h}\|_1 \\
&\lesssim \|\partial_t\vec{h}\|+\|\mathfrak{A}\vec{h}\|+\|\vec{h}\|\\
&=\|\partial_t\vec{h}\|+\|
\partial_t\vec{h}-\vec{f}\|+\|\vec{h}\| \\
&\lesssim \|\partial_t\vec{h}\|+\|\vec{h}\|+\|\vec{f}\|,
\end{align*}
thanks to Proposition 12. Suppose (\ref{H11}) is true for $n\leftleftarrows n-1$ such that
$n+2\leq\sigma$. Then, thanks to Proposition 12, we see
\begin{align*}
Z(\vec{h};n+1)&=Z(\partial_t\vec{h},n)+\|\vec{h}\|_{n+1} \\
&\lesssim \|\partial_t^{n+1}\vec{h}\|+\sum_{\nu\leq n-1}\|\partial_t^{\nu+1}\vec{h}\|+\|\vec{f}\|_{n-1} \\
&+\|\mathfrak{A}\vec{h}\|_n+\|\vec{h}\| \\
&\lesssim \|\partial_t^{n+1}\vec{h}\|+\sum_{\nu\leq n-1}\|\partial_t^{\nu+1}\vec{h}\|+\|\vec{f}\|_{n-1} \\
&+\|\partial_t\vec{h}-\vec{f}\|_n+\|\vec{h}\| \\
&\lesssim\|\partial_t^{n+1}\vec{h}\|+
\sum_{\nu\leq n}\|\partial_t^{\nu}\vec{h}\|+
\|\vec{f}\|_n,
\end{align*}
that is, (\ref{H11}) holds for $n$. Hence (\ref{H11}) holds for $n\leq\sigma-2$.

(\ref{H11}) implies
$$W(\vec{h};n+1)\lesssim \|\partial_t^{n+1}\vec{h}\|+W(\vec{h}; n)+ \|\vec{f}\|_n.$$
By induction on $n$, we have
\begin{equation}
W(\vec{h};n+1)\lesssim
\|\partial_t^{n+1}\vec{h}\|+\|\vec{h}\|+\|\vec{f}\|_n. \label{Z00}
\end{equation}

On the other hand, $\vec{H}=\partial_t^{n+1}\vec{h}$ satisfies
$$\frac{\partial\vec{H}}{\partial t}+
\mathfrak{A}\vec{H}=\vec{f}_{n+1}$$
with
$$ \vec{f}_{n+1}:=\partial_t^{n+1}\vec{f}-[\partial_t^{n+1},\mathfrak{A}]\vec{h}.$$
The energy estimate gives
\begin{equation}
\|\partial_t^{n+1}\vec{h}\|\lesssim
\|\partial_t^{n+1}\vec{h}|_{t=0}\|+
\int_0^t\|\vec{f}_{n+1}(t')\|dt'.\label{H13}
\end{equation}

Let us estimate $\|\partial_t^{n+1}\vec{h}|_{t=0}\|$. 

First we note
\begin{equation}
X_0(\vec{h};0,k)=\|\vec{h}|_{t=0}\|_k=\|\vec{0}\|_k=0
\quad\forall k.
\label{H14}
\end{equation}

Let $j+k=n$. Then
\begin{align*}
X_0(j+1,k)&=\|\partial_t^{j+1}\vec{h}|_{t=0}\|_k \\
&=\|-\partial_t^j\mathfrak{A}\vec{h}+\partial_t^j\vec{f}|_{t=0}\|k \\
&=\|
-\mathfrak{A}(\partial_t^j\vec{h})+
[\partial_t^j,\mathfrak{A}]\vec{h}+
\partial_t^j\vec{f}|_{t=0}\|_k \\
&\leq X_0(j,k+1)+W_0(\vec{h};n)+W_0(\vec{f};n),
\end{align*}
thanks to Propositions 13, 14. This and (\ref{H14}) imply
\begin{equation}
X_0(\vec{h};j,k)\lesssim
W_0(\vec{h};n)+W_0(\vec{f};n)
\end{equation}
for $j+k=n+1$. Hence
$$W_0(\vec{h}; n+1)\lesssim W_0(\vec{h}; n)+W_0(\vec{f}; n),$$
which implies
$$W_0(\vec{h};n+1)\lesssim W_0(\vec{f}; n),$$
since $W_0(\vec{h};0)=0$. Since
$\|\partial_t^{n+1}\vec{h}|_{t=0}\|\leq W_0(\vec{h}; n+1)$, (\ref{H13})
reads
\begin{equation}
\|\partial_t^{n+1}\vec{h}\|
\lesssim W_0(\vec{f};n)+
\int_0^t
\|\vec{f}_{n+1}(t')\|dt'.\label{Z01}
\end{equation}

On the other hand, 
\begin{equation}
\|\vec{f}_{n+1}(t)\|=
\|\partial_t^{n+1}\vec{f}-[\partial_t^{n+1},\mathfrak{A}]\vec{h}\| 
\lesssim\|\partial_t^{n+1}\vec{f}\|+|\vec{h}; t, n+1\|,\label{Z02}
\end{equation}
thanks to Proposition 14. Here we recall
$$
|\vec{h};\tau, n\|=\sup_{0\leq t\leq\tau}W(\vec{h};n).$$

Summing up (\ref{Z00}), (\ref{Z01}),(\ref{Z02}), we have
\begin{align}
|\vec{h};t, n+1\|&\lesssim 1+
W_0(\vec{f};n)+\|\vec{f}\|_n+ \nonumber\\
&+\int_0^t\|\partial_t^{n+1}\vec{f}(t')\|dt'+
\int_0^t|\vec{h};t', n+1\|dt' \label{Z03}
\end{align}
Here we have supposed
$$\|\vec{h}^{[\mu]}\|\lesssim\int_0^t\|\vec{f}^{[\mu]}\|
\lesssim\int_0^t\|\vec{g}^{[\mu]}\|+\|\vec{h}^{[1-\mu]}\|
\lesssim
\int_0^T\|\vec{g}\|\lesssim 1,$$
recalling
$$\vec{f}^{[\mu]}=\vec{g}^{[\mu]}+
\begin{bmatrix}
(-1)^{\mu}c_{11} & 0 \\
(-1)^{\mu}c_{21}\check{D}+(-1)^{\mu}c_{20} & c_{22}
\end{bmatrix}
\vec{h}^{[1-\mu]}.
$$
It follows from (\ref{Z03}) by the Gronwall's argument that
\begin{equation}
|\vec{h}^{[\mu]};t,n+1\|\lesssim
1+|\vec{f}^{[\mu]};t,n\|+
\int_0^t|\vec{f}^{[\mu]};t',n+1\|dt'.
\end{equation}

This implies
$$|\vec{h};T,n+1\|\lesssim
1+|\vec{g};T,n\|+
\int_0^T|\vec{g};t,n+1\|dt,$$
or, keeping in mind (\ref{X01}), we have

\begin{Proposition}
It holds that
$$|\|\vec{h}\||_{n+1}\lesssim 1+|\|\vec{g}\||_{n+1},$$
provided that $$s_N+2\leq\sigma, \quad n+2\leq \sigma,\quad
\|(1-x)^{N/2}\|_{\sigma}<\infty,\quad
|\|\vec{w}\||_{\sigma+3}\lesssim 1,\quad
|\|\vec{g}\||_1\lesssim 1.$$
\end{Proposition}

\section{Main result}

Now we are ready to apply the Nash-Moser(-Schwartz) theorem 
to our problem.\\

Take $2\nu=n+1, \sigma=2\nu+2=n+3$. Then
$n+2\leq \sigma$ is satisfied, and $s_N+2\leq \sigma
\Leftrightarrow s_N\leq 2\nu$. Proposition 15 reads
that
$$\|\vec{h}\|_{\nu}^{(2)}\lesssim1+\|\vec{g}\|_{\nu+1}^{(2)}, $$
provided that
$$(|\|\vec{w}\||_{\sigma+3}\lesssim)
\|\vec{w}\|_{\nu+3}^{(2)}\lesssim 1,\quad
\|\vec{g}\|_1^{(2)}\lesssim 1,\quad
\|(1-x)^{N/2}\|_{2\nu+2}<\infty.$$
(Recall (\ref{X01})(\ref{X02}).)
Since the Nash-Moser(-Schwartz) theorem requires that
$$\|\vec{w}\|_{E_j}\lesssim 1,\quad \|\vec{g}\|_{F_j}<\infty
\Rightarrow \|\vec{h}\|_{E_{j-1}}<\infty,$$
we should guarantee that
$\nu=\mathfrak{b}_E+\mathfrak{r}(j-1)$ satisfy $\nu+1\leq \mathfrak{b}_F+\mathfrak{r}j$ and
$\nu+3\leq \mathfrak{b}_E+\mathfrak{r}j$. That is, we require
\begin{equation}
\mathfrak{b}_E+\mathfrak{r}(j-1)+1\leq \mathfrak{b}_F+\mathfrak{r}j \label{H16}
\end{equation}
and
\begin{equation}
\mathfrak{b}_E+\mathfrak{r}(j-1)+3\leq \mathfrak{b}_E+\mathfrak{r}j.\label{H17}
\end{equation}
Now (\ref{H17}) is satisfied, if we take $\mathfrak{r}=3$ to fix the idea.
(\ref{H16}) reads
$$\mathfrak{b}_E-2\leq \mathfrak{b}_F$$
for $\mathfrak{r}=3$. Recall that (\ref{Sch2}) required
$$\mathfrak{b}_F\leq \mathfrak{b}_E-2.$$
So, (\ref{Sch2}) and (\ref{H16}) are satisfied if we take
$$\mathfrak{b}_F=\mathfrak{b}_E-2.$$

We should have 
$$s_N\leq 2\nu=2(\mathfrak{b}_E+\mathfrak{r}(j-1))$$
for $j=1,\cdots, 10$. It holds if $2\mathfrak{b}_E\geq s_N$.
If $N>6$, then $s_N\geq 4$, and $2\mathfrak{b}_E\geq s_N$ and
(\ref{Sch1}) are satisfied, if we take
$$\mathfrak{b}_E=s_N-2.$$
Summing up, we take
$$\mathfrak{b}_E=s_N-2,\quad \mathfrak{b}_F=s_N-4,\quad \mathfrak{r}=3, $$
provided that $N>6$.

Finally $\|(1-x)^{N/2}\|_{2\nu+2}<\infty$ should hold
for $\nu=\mathfrak{b}_E+\mathfrak{r}(J-1)=s_N+25$, where $J=10$. This means
\begin{equation}
\frac{3N}{2}>2s_N+52,\label{H18}
\end{equation}
which is equivalent to (\ref{Sch6}). This condition (\ref{H18}) is satisfied if $N>108$.

Thus we have

\begin{Theorem}
If $N>108$, the Nash-Moser(-Schwartz) theorem can be applied: If
$\|\mathfrak{F}(\vec{0})\|_{s_N+1}^{(2)}$ is sufficiently
small, there exists a solution $\vec{w}=(\tilde{y},\tilde{v})^T$ with
small $\|\vec{w}\|_{s_N+1}^{(2)}$. Note that
\begin{align*}
\mathfrak{F}(\vec{0})&=\begin{bmatrix}
c_1 \\
c_2
\end{bmatrix}, \\
c_1&=-\frac{\partial y^*}{\partial t}+J\Big(x,y^*, x\frac{\partial y^*}{\partial x}\Big)v^*, \\
c_2&=-\frac{\partial v^*}{\partial t}
-H_1\Big(x,y^*, x\frac{\partial y^*}{\partial x}, v^*\Big)\mathcal{L}y^*
-H_2\Big(x, y^*, x\frac{\partial y^*}{\partial x},
v^*, x\frac{\partial v^*}{\partial x}\Big)
\end{align*}
\end{Theorem}

\section{Remark for the application to the gaseous stars}

If we consider the gaseous stars governed by the Euler-Poisson
equations under the exact $\gamma$-law $P=A \rho^{\gamma}$,
the equilibrium governed by the Lane-Emden equation
has a finite radius if and only if
$\gamma>6/5$, that is, $N<12$. Therefore when $N$ is large,
the equilibriua cannot have finite radii. But the equation of state is not
exact $\gamma$-law, equilibria can have finite radii, even if
$\gamma-1$ near the vacuum
is very small. We shall show this.\\

Let us consider the Tolman-Oppenheimer-Volkoff equation
$$
\frac{dm}{dr}=4\pi r^2\rho, \qquad
\frac{du}{dr}=-\frac{G(m+4\pi r^2P/c^2)}{r^2(1-2Gm/c^2)},
$$
where $\rho$ and $P$ are give functions of $u$. Since we
define
$$u=\int_0^{\rho}\frac{dP}{\rho+P/c^2}$$
when the function  $\rho \mapsto P$ is given so that $ dP/d\rho=O(\rho^{\gamma-1})$ with $\gamma >1$, we should put
$$P=\int_0^u\varphi(u')e^{(u-u')/c^2}du' $$
when $\rho=\varphi(u)$ is given. Here $\varphi(u)$ is a positive smooth function of
$u>0$ such that 
$D\varphi(u)>0$ for $u>0$ and $\varphi(u)\rightarrow 0$ as $u\rightarrow +0$.
Suppose that there is a positive smooth function $\nu(u)$ of
$u>0$ such that $\nu(u)=\nu_0+[u]_1$ as $u\rightarrow +0$ with a constant
$\nu_0>1$ such that 
$$\rho=K^*\exp\Big[\int_{u*}^u\frac{\nu(u')}{u'}du'\Big].$$
It holds for $\rho=\varphi(u)$, when $\nu(u)=uD\varphi(u)/\varphi(u), K^*=\varphi(u^*)$, $u^*>0$ being arbitrary.
Note that if $$\nu(u)=\nu_0=\frac{N_0}{2}-1=\frac{1}{\gamma_0-1}=
\mbox{Const.}$$
for $0<u\ll 1$, then we have
$\rho\propto u^{\nu_0}$ and 
$P\propto \rho^{\gamma_0}(1+[\rho^{\gamma_0-1}]_1)$
as $\rho\rightarrow +0$.

Now let us consider $\rho=\rho_1(u)=u^{\nu_1}$ with a constant
$\nu_1=1/(\gamma_1-1)$ such that $1<\nu_1<3$ or $4/3<\gamma_1<2$.
Let $(m,u)=(m_1(r), u_1(r))$ be the solution such that
$(m, u)=(0,u_{1c})$ at $r=0$. Since $P\propto \rho^{\gamma_1}(1+O(\rho^{\gamma_1-1}))$ with $4/3<\gamma_1<2$,
\cite[Theorem 1]{TM98} says that this solution is short, that is,
there is a finite $r_{1+}$ such that 
$u_1(r) \searrow 0$ as $r\nearrow r_{1+}$. Hence we can find $r^*\in ]0,r_{1+}[$ such that $x_1(r^*)>1/G$, where $x=x_1(u)$ is defined by
$x=-m/ru$ along the solution $u=u_1(r)$. Of course we assume that $u_{1c}$ is so small that
$dP/d\rho < c^2$ for $0<u< 2u_{1c}$. Moreover we can assume that $r^*$ is
independent of large $c$, for,  as $c\rightarrow +\infty$, the Tolman-Oppenheimer-Volkoff equation
approaches the Lane-Emden equation
$$\frac{dm}{dr}=4\pi r^2u^{\nu_1},\qquad
\frac{du}{dr}=-\frac{Gm}{r^2}.
$$ 

Clearly there is a positive smooth function $\tilde{\nu}(u)$
of $u>0$ such that $\tilde{\nu}(u)=\nu_1$ for $u\geq u^*:=
u_1(r^*)$ and
$\tilde{\nu}(u)=\nu_0$ for $u\leq u^*/2$. 
Here we take $N_0$ arbitrarily
large and put $$\nu_0=\frac{N_0}{2}-1=\frac{1}{\gamma_0-1}.$$
Put
\begin{align*}
\tilde{\rho}(u)&=\rho_1(u^*)\exp\Big[\int_{u^*}^u
\frac{\tilde{\nu}(u')}{u'}du'\Big], \\
\tilde{P}(u)&=\int_0^u\tilde{\rho}(u')e^{(u-u')/c^2}du'.
\end{align*}
Let $(m,u)=(m_0(r), u_0(r))$ be the solution for $\rho=\tilde{\rho}(u), P=\tilde{P}(u)$ such that $(m,u)=(0,u_{1c})$. Since $\tilde{\rho}(u)=\rho_1(u)$ for $u\geq u^*$ and $\tilde{P}(u)=O(1),
d\tilde{P}/d\tilde{\rho}=O(1)$ as $c\rightarrow +\infty$,
we have $d\tilde{P}/d\tilde{\rho} < c^2$ and $|m_0(r)-m_1(r)|+|u_0(r)-
u_1(r)|=O(1/c^2)$ for $0\leq r\leq r^*$, so, $x_0(r^*)>1/G$
provided that $c$ is sufficiently large. 
In fact, as $c\rightarrow +\infty$, the Tolman-Oppenheimer-Volkoff equation
approaches the generalized Lane-Emden equation
$$\frac{dm}{dr}=4\pi r^2\tilde{\rho}(u),\qquad
\frac{du}{dr}=-\frac{Gm}{r^2}
$$
with error $O(1/c^2)$.
Then by \cite[p.61, Remark]{TM98}
we can claim that the solution $u_0(r)$ is short, that is, there is
a finite $r_{0+} <r^* \exp[1/(Gx_0(r^*)-1)]$ such that $u_0(r) \searrow 0$ as $r\nearrow r_{0+}$. Since
$P\propto \rho^{\gamma_0}(1+[\rho^{\gamma_0-1}]_1)$ with $N_0$ arbitrarily large, this is a desired 
example.\\

{\bf\large Acknowledgment}

A part of this work was done during the stay of the author at Nara Women's University at the opportunity of the Workshop `Free Boundary Problems in Fluid and Plasma Dynamics' in February, 2016. The stay was financially supported by JSPS KAKENHI Grant Number 15K04957.
The author would like to express his sincere thanks to the organizer  Professor Taku Yanagisawa (Nara WU ) for his kind hospitality, and  Professor Paolo Secchi (Brescia University) for his kind advices during discussions in the workshop.
Moreover, the author would like to express his sincere thanks to the 
anonymous
referee for pointing 
insufficiency
 of the original exposition. Thanks to his/her suggestions, the author could 
supplement
the manuscript. \\

{\bf\Large Appendix}\\

Let us describe the spectral properties of the linear
operator $\mathcal{L}$. We write 
$$
\mathcal{L}=-\frac{1}{b(x)}\frac{d}{dx}a(x)\frac{d}{dx}+L_0(x), 
$$
where
\begin{align*}
a(x)&=x^{5/2}(1-x)^{N/2}M(x), \\
b(x)&=x^{3/2}(1-x)^{N/2-1}M(x), \\
M(x)&=\exp\Big[\int_0^x\frac{L_1(x')}{x'(1-x')}dx'\Big]
\end{align*}
Note that $M(x)$ is a smooth function of $0\leq x <1$ and enjoys
$$
M(x)=\begin{cases}
M_0+[x]_1\quad\mbox{as}\quad x\rightarrow 0 \\
M_1+[1-x,(1-x)^{N/2}]_1\quad\mbox{as}\quad x\rightarrow 1
\end{cases}
$$
with positive constants $M_0, M_1$.

The Liouville transformation
\begin{align*}
\xi&=\int_{1/2}^x\frac{dx}{\sqrt{x(1-x)}}=\arcsin (2x-1),\\
y&=x^{-1}(1-x)^{-(N-1)/4}M(x)^{-1/2}\eta
\end{align*}
turns the equation
$$
\mathcal{L}y=\lambda y+f
$$
to the standard form
$$
-\frac{d^2\eta}{dx^2}+q\eta=\lambda \eta+\hat{f},
$$
where
$$
\hat{f}(x)=x^{-1/2}(1-x)^{-(N-3)/4}M(x)^{-1/2}f(x),
$$
and
$$q=L_0+\frac{1}{4}\frac{a}{b}
\Big(D\Big(\frac{Da}{a}+\frac{Db}{b}\Big)-\frac{1}{4}
\Big(\frac{Da}{a}+\frac{Db}{b}\Big)^2
+\frac{Da}{a}\Big(\frac{Da}{a}+\frac{Db}{b}\Big)\Big).
$$
Note that $x=0,1/2,1$ are mapped to $\xi=-\pi/2, 0, \pi/2$ and
$$x\sim
\begin{cases}
\displaystyle\frac{1}{4}\frac{1}{(\xi+\pi/2)^2} \quad\mbox{as}\quad x\rightarrow 0 \\
\\
\displaystyle\frac{1}{4}\frac{1}{(\pi/2-\xi)^2}\quad\mbox{as}\quad x\rightarrow 1.
\end{cases}
$$
We see
$$
q\sim
\begin{cases}
\displaystyle\frac{2}{(\xi+\pi/2)^2}\quad\mbox{as}\quad \xi\rightarrow -\pi/2 \\
\\
\displaystyle\frac{(N-1)(N-3)}{4}\frac{1}{(\pi/2-\xi)^2}
\quad\mbox{as}\quad \xi\rightarrow \pi/2.
\end{cases}
$$
Note that $(N-1)(N-3)/4 >3/4$ under the assumption {\bf(B0)}: $N\geq 5$.

Therefore we have \\

{\it 
The operator $\mathfrak{S}_0,\mathsf{D}(\mathfrak{S}_0)=
C_0^{\infty}(-\pi/2,\pi/2), $

\noindent $\mathfrak{S}_0\eta=-d^2\eta/dx^2+q\eta$ in
$L^2(-\pi/2,\pi/2)$ has the Friedrichs extension
$\mathfrak{S}$, a self-adjoint operator,
whose spectrum consists of simple eigenvalues
$\lambda_1<\cdots<\lambda_n<\lambda_{n+1}<\cdots\rightarrow
+\infty$. Thus the operator
$\mathfrak{T}_0, \mathsf{D}(\mathfrak{T}_0)=
C_0^{\infty}(0,1), \mathfrak{T}_0y=\mathcal{L}y$
in

\noindent $L^2(([0,1], x^{3/2}(1-x)^{N/2-1}dx)$ has the self-adjoint
extension $\mathfrak{T}$ with eigenvalues $(\lambda_n)_n$.} \\

Moreover we can claim\\

{\it If $\Phi(x)$ is an eigenfunction of the operator $\mathfrak{T}$,
we have
$$
\Phi(x)=
\begin{cases}
C_0+[x]_1\quad\mbox{as}\quad x\rightarrow 0 \\
C_1+[1-x, (1-x)^{N/2}]_1\quad\mbox{as}\quad x\rightarrow 1,
\end{cases}
$$
where $C_0, C_1$ are non-zero constants.}

\end{document}